\documentclass[journal]{IEEEtran}
\usepackage{comment}
\usepackage{cite}
\ifCLASSINFOpdf
  \usepackage[pdftex]{graphicx}
\else
\fi
\usepackage{amsmath}
\usepackage{amsfonts,amssymb,mathtools}
\allowdisplaybreaks
\newcommand*\Bell{\ensuremath{\boldsymbol\ell}}

\usepackage{algorithmic}
\usepackage{algorithm}

\usepackage{array}
\usepackage{booktabs}
\usepackage{multirow}

\ifCLASSOPTIONcompsoc
  \usepackage[caption=false,font=normalsize,labelfont=sf,textfont=sf]{subfig}
\else
  \usepackage[caption=false,font=footnotesize]{subfig}
\fi

\usepackage[colorlinks=true, bookmarks=true]{hyperref}
\AtBeginDocument{%
  \hypersetup{
    citecolor=red,
    linkcolor=blue,   
    urlcolor=blue}}

\begin{document}

\title{Reconfiguration and Real-Time Operation of Networked Microgrids Under Load Uncertainty}

\author{Hannah~Moring,
        Bala~Kameshwar~Poolla, 
        Harsha~Nagarajan,
        Johanna~L.~Mathieu, \\
        Andrey~Bernstein,
        and~David~M.~Fobes,
\thanks{This work was authored in part by the National Renewable Energy Laboratory, operated by Alliance for Sustainable Energy, LLC, for the U.S. Department of Energy (DOE) under Contract No. DE-AC36-08GO28308. Funding for the NREL \& LANL authors was provided by U.S. Department of Energy Office of Electricity's Microgrid Research and Development program under the ``Dynamic Microgrids for Large-Scale DER Integration and Electrification (DynaGrid)'' project. The views expressed in the article do not necessarily represent the views of the DOE or the U.S. Government. The U.S. Government retains and the publisher, by accepting the article for publication, acknowledges that the U.S. Government retains a nonexclusive, paid-up, irrevocable, worldwide license to publish or reproduce the published form of this work, or allow others to do so, for U.S. Government purposes.}
\thanks{H. Moring and J. L. Mathieu are with Electrical Engineering and Computer Science, University of Michigan, Ann Arbor, Michigan, USA (e-mail: \href{mailto:hmoring@umich.edu}{hmoring@umich.edu}, \href{mailto:jlmath@umich.edu}{jlmath@umich.edu}).}
\thanks{B. K. Poolla and A. Bernstein are with National Renewable Energy Laboratory, Golden, CO, USA (e-mail: \href{mailto:bpoolla@nrel.gov}{bpoolla@nrel.gov}, \href{mailto:Andrey.Bernstein@nrel.gov}{Andrey.Bernstein@nrel.gov}).}
\thanks{H. Nagarajan is with the Applied Mathematics \& Plasma Physics (T-5) group, Los Alamos National Laboratory (LANL), Los Alamos, NM, USA (e-mail: \href{mailto:harsha@lanl.gov}{harsha@lanl.gov}).}
\thanks{D. M. Fobes is with the Information Systems \& Modeling (A-1) group, LANL, Los Alamos, NM, USA (e-mail: \href{mailto:dfobes@lanl.gov}
{dfobes@lanl.gov}).}
}

\maketitle

\begin{abstract}
Distribution networks are increasingly exposed to threats such as extreme weather, aging infrastructure, and cyber risks--resulting in more frequent contingencies and outages, a trend likely to persist. Microgrids, particularly dynamic networked microgrids (DNMGs), offer a promising solution to mitigate the impacts of such contingencies and enhance resiliency. However, distribution networks present unique challenges due to their unbalanced nature and the inherent uncertainty in both loads and generation. This paper builds upon our prior work on the two-stage mixed-integer robust optimization problem for configuring DNMGs, improving the solve time and scalability. Furthermore, we present a model-free, real-time optimal power flow algorithm to manage DNMG operations in the time between reconfigurations. A case study on a realistic network based on part of the San Francisco Bay Area demonstrates the scalability of both approaches. The case study also illustrates the ability to maintain power flow feasibility as loads vary and operating conditions change when the methods are used in tandem. 
\end{abstract}

\begin{IEEEkeywords}
Dynamic Networked Microgrids, Uncertain Loads, Robust Optimization, Model-free Optimal Power Flow
\end{IEEEkeywords}

\section{Introduction}
More than a decade ago, the U.S. Department of Energy identified microgrids as a key building block in grids of the future~\cite{ton2012us}. Recent research has shown that microgrids can improve the resiliency and economic efficiency of distribution grids, especially during extreme weather events and unforeseen contingencies~\cite{Hamidieh2022Microgrids,Li2017Networked,barnes2019resilient}. Real-world implementations have shown that microgrids improve resilience during natural disasters, like Superstorm Sandy and Hurricane Maria~\cite{abbey2014powering,ferrari2023well}.

Beyond stand-alone microgrids, recent research has focused on what are known as nested microgrids, collaborative microgrids, or dynamic networked microgrids (DNMGs)~\cite{sharma2024engineering, Ma2018Real,Nassar2016Adaptive}. DNMGs are individual microgrids which can be connected to each other in order to facilitate resource sharing and form new microgrids with dynamic boundaries~\cite{zhu2021smart}. Compared to more traditional microgrids, DNMGs offer superior adaptability to volatile loading and uncertain generation availability~\cite{Bordbari2024Networked}. Their ability to reconfigure in response to contingencies, isolating impacted areas while maintaining operation in other parts of the network, can improve system resilience. Moreover, the ability to optimize over both connectivity and utilization of distributed generators (DGs) within DNMGs can reduce power losses and costs~\cite{Alam2019Networked}.

Despite the attention that DNMGs, or microgrids more generally, have received, there are still a number of open questions as to how to operate them~\cite{fobes2022optimal,Alam2019Networked}. The literature is missing a unified approach for addressing load uncertainty, contingency management, phase unbalance, network reconfiguration, inverter control modes, and real-time operation of DNMGs. Reconfiguration of distribution networks with various sources of uncertainty has been well studied~\cite{Zhou2022Three,Lee2015Robust,Mahdavi2023Robust}.  
However, these works do not consider the formation of microgrids or the need for load shedding. Distribution network reconfiguration with the ability to form microgrids is considered in~\cite{barani2018optimal,arefifar2012supply}. However, both works focus on the placement of switches and DGs rather than reconfiguration given fixed switch and DG locations. With the exception of~\cite{Zhou2022Three}, these papers and many others on distribution network reconfiguration~\cite{babaei2020distributionally,Arif2017Networked} consider only balanced networks despite the unbalanced nature of distribution networks. In~\cite{moring2024robust}, we addressed the problem of reconfiguration of DNMGs to maximize uncertain-load delivery in unbalanced multi-phase networks during a contingency using robust optimization and a cutting-plane method. A similar problem was approached using stochastic optimization, where load and solar generation scenarios represent uncertainty, in~\cite{Arif2017Networked}. However,~\cite{Arif2017Networked} considers balanced networks, assumes that individual loads can be shed, and does not address the role of grid-forming/grid-following inverter control modes.

Deciding network topology is only the first step in DNMG operation. Once the topology is chosen, i.e., where the microgrid boundaries lay is decided, operational decisions like power injections must be made. Existing literature on real-time operation of microgrids does not consider the complex reconfiguration enabled by DNMGs. For example, managing the energy exchanged between fixed-boundary microgrids with load uncertainty~\cite{Fathi2013Adaptive} or both load and generation uncertainty~\cite{Gao2018Decentralized, Nikmehr2015Optimal} given full network information has been considered. In~\cite{Gholami2019Proactive}, a two-stage robust formulation for scheduling DGs to mitigate load shedding during an islanding event is proposed. In~\cite{Feng2023Distributed}, a decentralized droop control scheme which relies on an augmented Newton-type power flow formulation is presented. These existing works assume a complete model of the distribution network, which may be unrealistic for a network composed of DNMGs with changing boundaries. The use of a model-free control strategy enables a distributed control approach that does not require sharing an updated network model with each DG every time the network is reconfigured.  

In this paper, we build on the work in~\cite{moring2024robust} to improve the inverter control mode criteria, solve time, and scalability.
We also incorporate a model-free, real-time optimal power flow (MFRT-OPF) algorithm that uses voltage measurements and measured changes in slack bus power injections to respond to load changes and manage voltage deviations via decentralized control actions that can be leveraged to operate the network between reconfigurations.

This paper contributes to the literature by (i) improving the solve time and scalability of the work in~\cite{moring2024robust}, (ii) introducing a novel maximum-phase connection model for grid-forming inverter eligibility of DGs that couples inverter-mode decisions with network topology via \textit{linear} one- and two-hop constraints constructed from already existing formulation variables, 
and (iii) integrating and adapting a previously-developed model-free control method to manage real-time load changes and voltage deviations within DNMGs that ensures network feasibility between reconfigurations given load uncertainty. We demonstrate the scalability of the proposed algorithms and quantify the reliability and resilience improvements they enable under load uncertainty on a realistic distribution network (based on part of the San Francisco Bay Area). 

The paper is organized as follows: Section~\ref{sec:RPPF} introduces the robust partitioning problem, with notation, network configuration constraints, and the single-stage formulation in Subsections~\ref{subsec:notation},~\ref{subsec:network_config}, and~\ref{subsec:single_stage}, respectively. Section~\ref{sec:cutting-plane} details the cutting-plane algorithm for solving the robust partitioning problem. The MFRT-OPF algorithm is described in Section~\ref{sec:rt_opf}. Section~\ref{sec:case_study} presents a case study using the two algorithms, and Section~\ref{sec:conclusion} provides the conclusions.

\section{Robust Partitioning Problem Formulation} \label{sec:RPPF}
This section introduces notation; network configuration constraints, including new reconfiguration constraints; and a modified single-stage robust optimization formulation from our prior work~\cite{moring2024robust}.

\subsection{Notation and Preliminaries}
\label{subsec:notation}
In this paper, vectors are indicated by bold typeface and matrices are indicated by blackboard bold typeface. The $\lvert \cdot \rvert$ operation represents cardinality when the input is a set, absolute value when the input is a number, and magnitude when the input is a complex number.

Consider a distribution network with a set of buses $\mathcal{N}$, a set of lines $\mathcal{L}$, a set of phases $\Phi$, and a set of switches $\mathcal{E}^\mathrm{sw}$. The network can be reconfigured by opening or closing switches. For generality, we assume the network has a meshed structure when all switches are closed. To align with typical protection coordination schemes~\cite{Lee2015Robust}, we assume that the network can only be operated in radial configurations. Within the network, there is a set of transformers $\mathcal{E}^x$ composed of a set of wye-connected transformers $\mathcal{E}^{x,Y}$ and a set of delta-connected transformers $\mathcal{E}^{x,\Delta}$, such that $\mathcal{E}^x=\mathcal{E}^{x,Y}\cup \mathcal{E}^{x,\Delta}$. Within the network, there is a set $\mathcal{G}$ of controllable DGs $g$. 
Let $s_g^\phi$ denote the apparent power injected by generator $g$ on phase $\phi$.

Additionally, there is a set of uncontrollable and uncertain loads $\mathcal{D}$. The apparent power demanded by load $d$ on phase $\phi$ is given by $s^\phi_d = s_d^{0,\phi} + u_{d}^{\phi}$, where $s_d^{0,\phi}$ is the nominal load power and $u_{d}^{\phi}$ is the load uncertainty defined by the set
\begin{multline}
     \mathcal{U}= \biggl\{ u_{d}^{\phi} \in \mathbb{R}^{|\mathcal{D} \times \Phi|} ~ \bigg| \biggr. ~
      \underline{s}^\phi_{d} - s_{d}^{0,\phi} \leqslant u_{d}^{\phi} \ \leqslant \  \overline{s}^\phi_{d} - s_{d}^{0,\phi}  \\
    ~ \forall \phi \in \Phi, \forall d \in \mathcal{D} \biggr\},
    \label{eq:uncertaintySet}
\end{multline}
where $\overline{s}^\phi_{d}$ and $\underline{s}^\phi_{d}$ represent the upper and lower bounds of uncertain load $d$ on phase $\phi$.

The formulation presented in this paper uses the \textsc{LinDist3Flow} model to approximate the AC power flow equations while capturing the unbalanced, three-phase power flow typical of distribution grids~\cite{Gan2014Convex,Fobes2020PMD}. We adopt this approximation because its linear—and therefore convex—structure converts the inner dispatch problem in our robust max–min formulation (Section \ref{sec:cutting-plane}) into a continuous convex program, enabling the use of strong duality and ensuring computational tractability. Section \ref{subsec:RPOP-results} confirms that, for the network studied, this approximation remains AC-feasible to high accuracy, validating its use in the proposed optimization scheme. 

In the \textsc{LinDist3Flow} model, the equations are written in terms of squared voltage and
current magnitudes.
Let $w^\phi_i$ represent the squared magnitude of the voltage on phase $\phi$ at bus $i$ and let $l^\phi_{ij}$ represent the squared magnitude of the current flowing on phase $\phi$ of the line connecting buses $i$ and $j$. Let $z^{\phi\psi}_{ij} = r^{\phi\psi}_{ij} + jx^{\phi\psi}_{ij}$ represent the self and mutual impedances on the line connecting buses $i$ and $j$ between phases $\phi$ and $\psi$. When $\phi = \psi$, $z$ represents self impedance, otherwise it represents mutual impedance. The complex power flowing on line $ij$ on phase $\phi$ is given by $s^\phi_{ij} = p^\phi_{ij} + jq^\phi_{ij}$, where $p$ and $q$ represent active and reactive power, respectively.

\begin{figure}
    \centering
    \includegraphics[width=1\columnwidth]{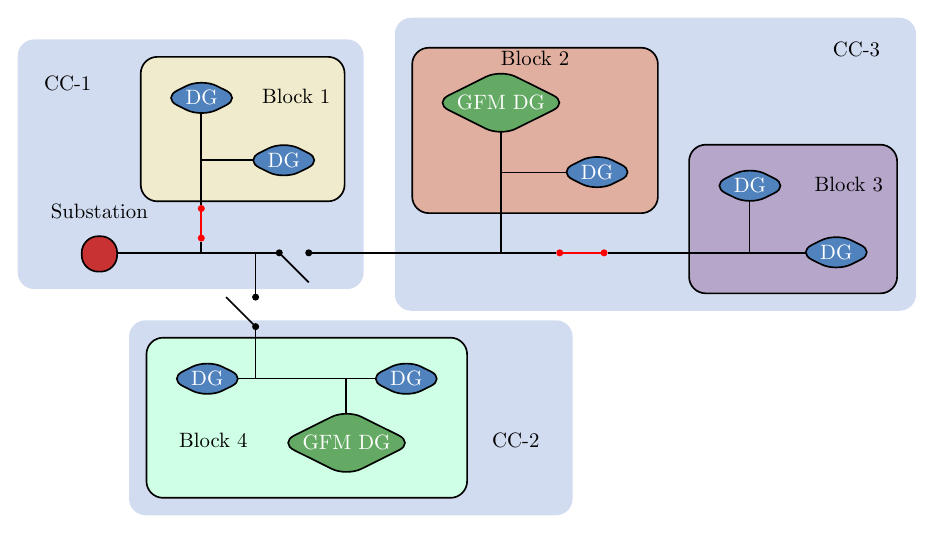}
    \caption{Dynamic reconfiguration of networked microgrids via switches. Blocks are connected by switches to form connected components (CCs). Each energized CC, which functions as an independent microgrid, contains uncertain loads and DGs, with one DG serving as a voltage source (GFM-DG).}
    \label{fig:load_block}
    \vspace{-.5cm}
\end{figure}

A connected component (CC) is a set of buses that are electrically connected, including single isolated buses~\cite{Lei2020Radiality}.
We consider a microgrid as a CC containing at least one load and at least one voltage source. We define $\mathcal{B}$ as the set of CCs that exist when every switch in the network is open. The CCs within this set will be referred to as blocks. When a switch between two blocks is closed, the two blocks form one CC, as shown in Fig.~\ref{fig:load_block}. Thus, the set of blocks in the network remains static, while the CCs change as switches open or close.

\subsection{Network Configuration Constraints}
\label{subsec:network_config}
Two general requirements govern network configuration: First, every energized CC must be a spanning tree, devoid of loops to satisfy typical protection schemes. Second, each energized CC must contain at least one voltage source.

\subsubsection{Radiality Constraints}
\label{subsubsec:radial}
A radial topology is enforced in each CC using the directed multi-commodity flow-based model of spanning tree constraints given in~\cite{Lei2020Radiality}. For brevity, we omit the constraints here and refer the reader to~\cite{Lei2020Radiality} for details. In the remainder of this paper, this model will be referred to as the radiality constraints.

\subsubsection{Grid-Forming DG Constraints}
An energized CC must contain a voltage source that is connected to each active phase. Distribution network DGs are typically inverter-based resources, which can operate as either grid-forming or grid-following. A grid-forming DG (GFM-DG) controls the AC-side voltage, acting as a voltage source for the network, while a grid-following DG controls the AC-side current and ``follows the phase angle of the existing grid voltage''~\cite{Li2022Revisiting}. Uncoordinated operation of two or more GFM-DGs in close proximity can lead to independent control of voltage and frequency, potentially leading to instability~\cite{Sadeque2021Multiple}. Control strategies exist for facilitating multiple GFM-DGs in a microgrid~\cite{Sharma2022Synchronization, Watson2021Scalable}, but may be incompatible with older, existing equipment~\cite{Han2016Review}. Consistent with the adopted model-free control strategy, each CC is restricted to include a single GFM-DG; for a generalized (multi-GFM) formulation, see~\cite{moring2024robust}. 

Let $z_l^\mathrm{bl}$, $z_{ij}^\mathrm{sw}$, and $z_g^\mathrm{inv}$ indicate whether block $l$ is energized, the switch between blocks $i$ and $j$ is closed, and DG $g$ is grid-forming, respectively. Each variable equals 1 if active (energized block, closed switch, or grid-forming DG) and 0 otherwise.
Let $\mathcal{G}_l$ and $\mathcal{E}^\mathrm{sw}_l$ be the set of DGs and switches, respectively, connected to block $l$. The relationship between the block states $z_l^\mathrm{bl}$ and the DG states $z_g^\mathrm{inv}$ is defined by
\begin{align}
& |z_i^\mathrm{bl} - z_j^\mathrm{bl}| \leqslant (1-z^\mathrm{sw}_{ij}),  \qquad \forall ij \in \mathcal{E}^\mathrm{sw} \label{eq:ConnectedBlocks} \\
& z^\mathrm{bl}_l \underline{s_g} \leqslant s_g^\phi \leqslant z^\mathrm{bl}_l \overline{s_g},  \qquad \forall g \in \mathcal{G}_l, \forall l \in \mathcal{B} \label{eq:BlockGenLimits} \\
& z^\mathrm{bl}_l - \sum_{ij \in \mathcal{E}^\mathrm{sw}_l} z^\mathrm{sw}_{ij} \leqslant \sum_{i \in \mathcal{G}_l} z^\mathrm{inv}_i \leqslant \mathrlap{z^\mathrm{bl}_l,}  \qquad \forall l \in \mathcal{B} \label{eq:GenPerBlock}
\end{align}
where~\eqref{eq:ConnectedBlocks} enforces that blocks connected by a closed switch must be either both energized or both de-energized; \eqref{eq:BlockGenLimits} enforces generation power limits for DGs within an energized block and prevents a DG within a de-energized block from producing power; and \eqref{eq:GenPerBlock} states that a block must contain exactly one GFM-DG if it is energized and not connected to another block. 

\paragraph{Coloring Scheme} To ensure each CC contains a GFM-DG and to ensure consistency across networked blocks, i.e., all blocks within a CC share the same
grid-forming inverter, we apply a coloring scheme in conjunction with a multi-commodity flow model, as implemented in PowerModelsONM.jl~\cite{fobes2022optimal}. In this coloring scheme, closed switches are colored, or labeled, with the block in their CC that contains the GFM-DG. The mathematical details of the coloring scheme and multi-commodity flow model are given in our prior work~\cite{moring2024robust} (see constraints (5)-(15)) and are omitted here for brevity. In the remainder of this paper, these constraints will be referred to as the coloring scheme constraints. 

\paragraph{Maximum Phase Connections} If the network contains blocks with different phase configurations, i.e., some are 3-phase while others are 2- or single-phase, additional constraints are needed to ensure that the DG selected as GFM is connected to every active phase in the CC. Each block $l \in \mathcal{B}$ is assigned a maximum phase parameter $\Phi_l^\text{max} \in \{1, 2, 3\}$, representing the highest number of phases among all buses in block $l$. Formally, 
\begin{equation*}
    \Phi_l^\text{max} = \max_{i \in \mathcal{N}_l} |\Phi_i |, \quad \forall l \in \mathcal{B}
\end{equation*}
where $|\Phi_i|$ indicates the number of phases present at bus $i$ and $\mathcal{N}_l$ is the set of buses in block $l$. 
Within each block $l$, only DGs connected to all $\Phi_l^{\max}$ phases are eligible to operate as GFM. Let $\mathcal{D}_l$ denote the set of DGs located in block $l$, and define $\mathcal{G}_l^{\Phi_l^{\max}} := \{\, g \in \mathcal{D}_l : g \text{ is connected to every phase present in } l \,\}$. We enforce this eligibility by requiring DGs outside this set to be grid-following:
\begin{equation}
z_g^{\text{inv}} = 0, \quad \forall g \notin \mathcal{G}_l^{\Phi_l^\text{max}}, \ \forall l \in \mathcal{B}.
\label{eq:GF_has_max_ph}
\end{equation} 

Because CCs can contain multiple blocks connected by switches, it must also be verified that any GFM-DG is connected to the maximum number of phases active within the blocks it is connected to. 
To formalize the resulting topological constraints, consider the following ``one-hop'' and ``two-hop'' constraints
\begin{align}
    & z_g^\text{inv} \leqslant 1 - z_{lj}^\text{sw}, \quad \forall g \in \mathcal{D}_l, \forall lj \in \mathcal{E}^\mathrm{sw,\underline{\Phi}}_{l}, \forall l \in \mathcal{B}  \label{eq:one-hop} \\
    \begin{split}
        z_g^\text{inv} \leqslant 2 - z_{lj}^\text{sw} -  z_{jk}^\text{sw}, \quad \forall g \in \mathcal{D}_l, \forall lj \in \mathcal{E}^\mathrm{sw}_l\setminus \mathcal{E}^\mathrm{sw,\underline{\Phi}}_{l},\\
         jk \in \mathcal{E}^\mathrm{sw,\underline{\Phi}}_{j}, \forall l \in \mathcal{B}
    \end{split} \label{eq:two-hop}
\end{align}
where $\mathcal{E}^\mathrm{sw,\underline{\Phi}}_{l}$ represents the set of switches that connect block $l$ to a neighboring block $j$ such that $\Phi_l^\text{max} < \Phi_j^\text{max}$. In words, \eqref{eq:one-hop}–\eqref{eq:two-hop} permit a DG in block $l$ to be GFM only if, under the current switch status, there is no block with $\Phi^\text{max}>\Phi_l^\text{max}$ that is either adjacent to $l$ or adjacent to a neighbor of $l$; otherwise the DG must be grid-following.

For chains of length greater than two consisting solely of blocks with $\Phi^\mathrm{max}<3$, the constraints \eqref{eq:one-hop}–\eqref{eq:two-hop} remain necessary but may not be sufficient; additional multi-hop analogs can be introduced by conditioning GFM feasibility on the closed status of longer paths to any block whose $\Phi^\mathrm{max}$ exceeds that of the DG’s home block. In summary, the proposed formulation generalizes the single-phase assumptions in~\cite{moring2024robust} to multi-phase networks by coupling inverter-mode decisions with the topology induced by switch statuses and the spatial variation of phase availability. 

\subsection{Single-Stage Robust Partitioning and Operation Problem}
\label{subsec:single_stage}
This section details the single-stage Robust Partitioning and Operation Problem (RPOP) for distribution grids composed of DNMGs. Its optimal solution comprises network partitions and generator set-points that minimize the costs of load shed and power generation while ensuring feasibility and robustness against all realizations of uncertain loads. The formulation is 
\begin{subequations} 
\label{eq:FullForm}
\begin{align}
     \min~ & \sum_{l\in\mathcal{B}}  \beta^\mathrm{B}_l(1-z^\mathrm{bl}_l)+ \sum_{\substack{g \in G}} \sum_{\substack{\phi \in \Phi_g}}c_{1,g}s^\phi_{g} + c_{0,g}z^\mathrm{bl}_{g}\label{obj} \\
     \mathrm{s.t.}~ & \text{Radiality~Constraints}, \text{Coloring Scheme~Constraints}, \nonumber \\ & \eqref{eq:ConnectedBlocks}-\eqref{eq:two-hop} \nonumber \\
    \begin{split}
         z^\mathrm{bl}_l\underline{v_{i}}^2 \leqslant w_{i,\phi} \leqslant z^\mathrm{bl}_l\overline{v_{i}}^2, ~\forall \phi \in \Phi_i, \\
         \forall i \in \mathcal{N}_l, \forall l \in \mathcal{B} \label{eq:voltLim}
     \end{split}  \\
     \begin{split}
         \mathbf{w}_i = \mathbf{w}_j - \mathbb{M}^P_{ij}\mathbf{p}_{ij} - \mathbb{M}^Q_{ij}\mathbf{q}_{ij}, \\
         \forall i,j \in \mathcal{N}, ~\forall ij \in \mathcal{L} \label{eq:VoltBal}
     \end{split}  \\
     & \lvert s^\phi_{ij} \rvert \leqslant \overline{s_{ij}^{\phi}}, ~\quad \forall \phi \in \Phi_{ij}, \forall ij \in \mathcal{L} \cup \mathcal{E}^{x} \label{eq:PowerLim}\\
     & \lvert s^\phi_{ij} \rvert \leqslant z^\mathrm{sw}_{ij}\overline{s_{ij}^{\phi}}, ~\quad \forall \phi \in \Phi_{ij}, \forall ij \in \mathcal{E}^\mathrm{sw} \label{eq:SWPowerLim}\\
     & w^\phi_i = n_{ij}^2 w^\phi_j, ~\quad \forall \phi \in \Phi_{ij}, \forall ij \in \mathcal{E}^{x,Y} \label{eq:XfmrYV}\\
     \begin{split}
         & 3(w_{i,\phi}+w_{i,\psi})= 2(n_{ij})^2w_{j,\phi}, \\
         &\quad \quad \forall (\phi,\psi) \in \{(a,b),(b,c),(c,a)\}, ~\forall ij \in \mathcal{E}^{x,\Delta} \label{eq:XfmrDeltaV}
     \end{split}\\
     & s^\phi_{ij} = s^\phi_{ji}, ~\quad \forall \phi \in \Phi_{ij},\forall ij \in \mathcal{E}^{x,Y} \label{eq:XfmrY} \\
    \begin{split}
        2p_{ij,\phi}=-(p_{ji,\phi}+p_{ji,\psi})+(q_{ji,\psi}-q_{ji,\phi})/\sqrt{3}, \\ ~\forall (\phi,\psi) \in \{(a,c),(b,a),(c,b)\}, ~\forall ij \in \mathcal{E}^{x,\Delta} \label{eq:XfmrDeltaP}
    \end{split}\\
    \begin{split}
        2q_{ij,\phi}=(p_{ji,\phi}-p_{ji,\psi})/\sqrt{3}-(q_{ji,\psi}+q_{ji,\phi}), \\
        ~\forall (\phi,\psi) \in \{(a,c),(b,a),(c,b)\}, ~\forall ij \in \mathcal{E}^{x,\Delta} \label{eq:XfmrDeltaQ}
    \end{split} \\
    & z^\mathrm{bl}_i\underline{s^\phi_{g}} \leqslant s^{\phi}_{g} \leqslant z^\mathrm{bl}_i\overline{s^\phi_{g}},~ \forall \phi \in \Phi_g, \forall g \in i, \forall i \in \mathcal{B} \label{eq:genLims} \\
    \begin{split}
         \sum_{ij \in \mathcal{E}_i} s^\phi_{ij} = \sum_{g \in i} s^\phi_{g}
         - z^\mathrm{bl}_k \sum_{d \in \mathcal{D}_i}s^\phi_d 
         - \sum_{c \in \Phi_i} y^\phi_c w^\phi_i, \\ ~\forall \phi \in \Phi_i,\forall i \in l,~\forall l \in \mathcal{B} \label{eq:PowerBalance}
     \end{split}
\end{align}
\end{subequations}
where $z^\mathrm{bl}_{g}$ represents the block state of the block containing generator $g$, and $\mathcal{E} = \mathcal{L}\cup\mathcal{E}^\mathrm{sw}\cup\mathcal{E}^{x}$ is the set of all lines, switches, and transformers.
The objective function minimizes the cost of de-energizing blocks and generation, where $\beta^\mathrm{B}_l$ is a weighting parameter for block $l$ indicating priority, and $c_{1,g}, c_{0,g}$ are parameters of the linear cost function of generator $g$. The voltage magnitude limits at each bus are enforced via~\eqref{eq:voltLim}. The linearized voltage drop between a downstream bus $i$ and its upstream bus $j$ is modeled by~\eqref{eq:VoltBal}. In this constraint, variables in bold denote three-phase vector quantities, while $\mathbb{M}^P_{ij}$ and $\mathbb{M}^Q_{ij}$ are complex-valued matrices that characterize the sensitivity of voltage drops to active and reactive power flows, respectively. These matrices are derived in~\cite{Arnold2016Optimal} (see eqns. (20) and (21) therein for the linearized three-phase voltage drop model).

Constraints~\eqref{eq:PowerLim} and~\eqref{eq:SWPowerLim} 
define the power flow limits for each line, transformer, and switch. 
Constraints~\eqref{eq:XfmrYV} and~\eqref{eq:XfmrDeltaV} define the voltage transformations across each wye-connected and delta-connected transformer, where $n_{ij}$ is the tap ratio of transformer $ij$. Constraints~\eqref{eq:XfmrY},~\eqref{eq:XfmrDeltaP} and~\eqref{eq:XfmrDeltaQ} define the relationship between the directions of power flow on each wye-connected and delta-connected transformer. Lastly,~\eqref{eq:PowerBalance} represents linearized power balance.

\section{A Cutting-Plane Algorithm to solve the RPOP} \label{sec:cutting-plane}
It is intractable to solve the RPOP as given in Section~\ref{subsec:single_stage} using off-the-shelf solvers. We now present a two-stage version of the RPOP, which can be solved to global optimality using an iterative cutting-plane algorithm. 

\subsection{Two-stage Robust Partitioning and Operation Problem}
\label{subsec:2stage_RPOP}
The RPOP in~\eqref{eq:FullForm} has an equivalent two-stage formulation composed of a master problem and a set of subproblems~\cite{Geoffrion1972Generalized}. The master (first-stage) decision variables include power injections, $\mathbf{s}^\phi_g~ \forall \phi \in \Phi,~\forall g \in \mathcal{G}$, switch configurations $z^\mathrm{sw}_{ij}~ \forall ij \in \mathcal{E}^\mathrm{sw}$, generator operating states $z^\mathrm{inv}_{g}~ \forall g \in \mathcal{G}$ as either grid-forming or grid-following, block states $z^\mathrm{bl}_l~ \forall l\in \mathcal{B}$, and additional configuration and power flow variables.

Let $\mathbf{x}$ represent a vector of variables $\mathbf{s}^\phi_g,z^\mathrm{sw}_{ij},z^\mathrm{inv}_{g},z^\mathrm{bl}_l$ including every generator, switch, and block in the network. Let $\mathbf{s}_{d}^*$ be the worst-case uncertain load realization.
The master problem, which relaxes the single-stage problem by ignoring load uncertainty, is
\begin{subequations} 
\begin{align}
    \min~ & \sum_{l\in\mathcal{B}}  \beta^\mathrm{B}_l(1-z^\mathrm{bl}_l)+ \sum_{\substack{g \in G}} \sum_{\substack{\phi \in \Phi_g}}c_{1,g}s^\phi_{g} + c_{0,g}z^\mathrm{bl}_{g} + \theta \nonumber \tag{$M$}\label{eq:master}\\
    \mathrm{s.t.}~ & \text{Radiality~Constraints}, \text{Coloring Scheme~Constraints}, \nonumber \\ & \eqref{eq:ConnectedBlocks}-\eqref{eq:two-hop} \nonumber \\
    & V_2(\mathbf{x}^*_k,\mathbf{s}_{d,k}^*) + \mathbf{\pi}_k^{\intercal} A(\mathbf{x}-\mathbf{x}^*_k) \leqslant \theta,
        \forall k=1,2,... \label{eq:Cuts} 
\end{align} 
\end{subequations} 
where subscript $k$ denotes the $k^\text{th}$ iteration; $\mathbf{x}^*_k$ is the optimal solution of~\eqref{eq:master}; $V_2(\mathbf{x}^*_k, \mathbf{s}_{d,k}^*)$ is the subproblem objective at iteration $k$; $\boldsymbol{\pi}$ is the vector of dual variables for the subproblem constraints involving master variables; and $A$ is coefficient matrix from~\eqref{eq:subproblem} corresponding to the master variables. Our previous work's master problem omitted constraints~\eqref{eq:voltLim}–\eqref{eq:PowerBalance}. By including the \textsc{LinDist3Flow} equations for a single load scenario, the initial master solution starts with a much tighter lower bound, which significantly reduces the number of iterations needed for convergence.

In the subproblem (second-stage), the integer variables and generation power set-points from the master problem are fixed. After the uncertainty is realized, adjustments to the generation power set-points, $o^{+}$ and $o^{-}$, are permitted but bounded by generator ramping capabilities and capacity constraints. The subproblem formulation is
\begin{subequations} 
\label{eq:subprob}
\vspace{-0.05in}
\begin{align}
    \begin{split} &\min~ \beta^\mathrm{S} \left( \sum_{i \in \mathcal{N}}\sum_{\phi \in \Phi_i}\left(h^{+,\phi}_i + h^{-,\phi}_i \right) \right)
    \\ & \qquad \qquad \qquad + \sum_{g \in G}\sum_{\phi \in \Phi_g} c_{1,g} \left(o^{+,\phi}_{g} - o^{-,\phi}_{g}\right) 
    \end{split} \nonumber \tag{$S_1$}\label{eq:subproblem}\\
    & \mathrm{s.t.}~ \quad \eqref{eq:voltLim}-\eqref{eq:XfmrDeltaQ}, \nonumber \\
    \begin{split}
         \sum_{ij \in \mathcal{E}_i} s^\phi_{ij} + h^{+,\phi}_i -  h^{-,\phi}_i = \sum_{g \in i} \left(s^{\phi*}_{g} + o^{+,\phi}_g - o^{-,\phi}_g \right)
         \\ - z^\mathrm{bl*}_k \sum_{d \in \mathcal{D}_i}s^\phi_d 
         - \sum_{c \in \Phi_i} y^\phi_c w^\phi_i, ~\forall i \in k,~\forall k \in \mathcal{B}, \label{eq:SubPowerBalance}
     \end{split} \\
     & o^{+,\phi}_{g} \leqslant z^\mathrm{bl*}_i\overline{s^\phi_{g}}  - s^{\phi*}_{g}, ~ \forall \phi \in \Phi_g, \forall g \in i, \forall i \in \mathcal{B},  \label{eq:SubModUpLims}\\
     & o^{-,\phi}_{g} \leqslant s_{g}^{\phi*} - z^\mathrm{bl*}_i\underline{s^\phi_{g}}, ~ \forall \phi \in \Phi_g,~\forall g \in i,~ \forall i \in \mathcal{B}, \label{eq:SubModDownLims}\\
     & 0 \leqslant \ o^{+,\phi}_{g}, \ o^{-,\phi}_{g} \ \leqslant \overline{o^\phi_{g}},~ \forall \phi \in \Phi_g,~\forall g \in i,~ \forall i \in \mathcal{B}, \label{eq:Sub+RampLims}  \\
     & h^{+,\phi}_i, h^{-,\phi}_i \geqslant 0 , ~ \forall \phi \in \Phi_i,~ \forall i \in \mathcal{N}. \label{eq:SubSlackPos}
\end{align}
\end{subequations}
where $\beta^\mathrm{S}$ is a weighting parameter and $h_i^{+,\phi}, h_i^{-,\phi}$ are slack variables indicating a power balance violation on phase $\phi$ at bus $i$. The realized load and injection adjustments are incorporated into~\eqref{eq:SubPowerBalance}. If any  slack variables are non-zero at optimality, there exists no feasible $o^{+},o^{-}$ that satisfy the power flow equations for the candidate solution $\mathbf{x}^*$ and uncertain load realization $\mathbf{s}_d$. In such cases, a feasibility cut \eqref{eq:Cuts} is generated to remove $\mathbf{x}^*$ from the feasible region of subsequent iterations of the master problem.
Constraints~\eqref{eq:SubModUpLims}-\eqref{eq:Sub+RampLims} set limits on generator set-point modifications based on ramping and capacity limits.

A decomposition algorithm cannot be directly implemented due to infinitely many subproblems arising from uncertain load realizations. Instead, we identify and solve only the subproblem with the most violated power balance constraints within the given uncertainty set. This entails solving a bilevel $\max-\min$ problem of the form
\begin{subequations}  
\begin{align}
    \begin{split} \max_{\mathcal{U}} & \min~ \beta^\mathrm{S} \left( \sum_{i \in \mathcal{N}}\sum_{\phi \in \Phi_i}\left(h^{+,\phi}_i + h^{-,\phi}_i \right) \right)
    \\ & \qquad \qquad \qquad + \sum_{g \in G}\sum_{\phi \in \Phi_g} c_{1,g} \left(o^{+,\phi}_{g} - o^{-,\phi}_{g}\right) 
    \end{split} \nonumber \tag{$S_2$}\label{eq:MaxMin} \\ 
    & \mathrm{s.t.} \quad \eqref{eq:voltLim}-\eqref{eq:XfmrDeltaQ}, \eqref{eq:SubPowerBalance}-\eqref{eq:SubSlackPos} \nonumber
\end{align}
\end{subequations} 

We know from convex optimization theory that the optimal solution to the $\max-\min$ problem~\eqref{eq:MaxMin} lies at one of the extreme points (minimum or maximum values) of the uncertainty set, since the inner minimization is a convex function of the uncertain load~\cite{Bertsekas2009Convex}. To improve the computational tractability of our prior work~\cite{moring2024robust}, we adopt a simplified uncertainty model by assuming that loads vary uniformly within predefined ``clusters''--i.e., all loads in a cluster simultaneously take on either their minimum or maximum values. This assumption is reasonable, as it captures correlated load behavior often observed in practice, such as geographically grouped loads responding similarly to weather, pricing signals, or time-of-use effects. By limiting the uncertainty set to a smaller set of scenarios, defined by the extreme points of these clusters, we reduce its dimensionality and significantly lower the number of subproblem evaluations required in the cutting-plane algorithm \ref{alg:cutting_plane}. Leveraging this, we replace \(\mathcal{U}\) in~\eqref{eq:MaxMin} with a finite-dimensional uncertainty set
\begin{multline} \label{eq:extreme_uncertaintySet}
     \widehat{\mathcal{U}}= \biggl\{ u_{d}^{\phi} \in \mathbb{R}^{|\mathcal{D} \times \Phi|} ~ \bigg| \biggr. ~
      u_{d}^{\phi} = \zeta_\gamma^{+}(\overline{s}^\phi_{d}) + \zeta_\gamma^{-}(\underline{s}^\phi_{d}) - s_{d}^{0,\phi},      \\
      \hspace{0.8cm} \zeta_\gamma^{+}+\zeta_\gamma^{-}=1, \ \zeta_\gamma^{+},\zeta_\gamma^{-} \in \{0,1\}, \\
    ~ \forall \phi \in \Phi, \forall d \in \mathcal{D}_\gamma, \forall \gamma \in \Gamma \biggr\}.
\end{multline}
where $\Gamma$ is a set of load clusters.
 Thus, the reduced finite-dimensional~\eqref{eq:MaxMin} problem can be replaced by solving the modified subproblem for every possible combination of $\zeta_d^+,\zeta_d^-$, selecting the solution that yields the largest objective, i.e., 
\begin{subequations} 
\label{eq:ExtremeSub}
\begin{align}
    \begin{split} &\min~ \beta^\mathrm{S} \left( \sum_{i \in \mathcal{N}}\sum_{\phi \in \Phi_i}\left(h^{+,\phi}_i + h^{-,\phi}_i \right) \right)
    \\ & \qquad \qquad \qquad + \sum_{g \in G}\sum_{\phi \in \Phi_g} c_{1,g} \left(o^{+,\phi}_{g} - o^{-,\phi}_{g}\right) 
    \end{split} \nonumber \tag{$E$}\label{eq:enumerate}\\  \\
    & \mathrm{s.t.} \quad \eqref{eq:voltLim}-\eqref{eq:XfmrDeltaQ}, \eqref{eq:SubModUpLims}-\eqref{eq:SubSlackPos}, \nonumber \\
    \begin{split}
         \sum_{ij \in \mathcal{E}_i} s^\phi_{ij}  + h^{+,\phi}_i -  h^{-,\phi}_i = \sum_{g \in i} \left(s^{\phi*}_{g} + o^{+,\phi}_g - o^{-,\phi}_g \right)
         \\ - z^\mathrm{bl*}_k \sum_{d \in \mathcal{D}_i} \big( s^{0,\phi}_d 
          + \zeta^+_\gamma\overline{s}^{\phi}_d - \zeta^-_\gamma\underline{s}^{\phi}_d \big) - \sum_{c \in \Phi_i} y^\phi_c w^\phi_i, \\ ~\forall i \in k,~\forall k \in \mathcal{B}_\gamma, \forall \gamma \in \Gamma, \nonumber
     \end{split} \\
     & \zeta^+_\gamma,\zeta^-_\gamma \in \{0,1\}, \forall \gamma \in \Gamma. \nonumber
\end{align}
\end{subequations}
Note that for simplicity we assume that all loads in a block are in the same load cluster. This approach still requires solving $2^{|\Gamma|}$  number of quadratically-constrained convex programs~\eqref{eq:enumerate}, and may not scale efficiently as the number of uncertain load clusters increases. 

\subsection{A Cutting-plane Algorithm}
We propose a cutting-plane algorithm, outlined in Algorithm~\ref{alg:cutting_plane}, to solve the two-stage RPOP. After initialization, the master problem~\eqref{eq:master} is solved for the network configuration and generator set-points (master solution, $\mathbf{x}^*$) in Step~\ref{alg:solve_M} of Algorithm~\ref{alg:cutting_plane}. Given $\mathbf{x}^*$, \eqref{eq:MaxMin} is solved for every scenario in $\widehat{\mathcal{U}}$ to obtain the worst-case load realization in Step~\ref{alg:solve_maxmin}. Finally, the subproblem~\eqref{eq:subproblem} is solved in Step~\ref{alg:solve_S} to check the feasibility of the master solution given the worst-case load realization. If $\mathbf{x}^*$ is feasible, then it is an optimal solution to the RPOP and the algorithm terminates. Otherwise, if there is non-zero slack $h^{+,\phi},  h^{-,\phi}$, a sub-gradient cut is added to~\eqref{eq:master} in Step~\ref{alg:add_cut} and the process is repeated.
\begin{algorithm}[ht]
  \caption{: Cutting-plane algorithm for two-stage RPOP}
  \label{alg:cutting_plane}
\begin{algorithmic}[1]
\STATE Initialize: $k \leftarrow 1$, $h^{+,\phi}_i \leftarrow \infty$, $h^{-,\phi}_i \leftarrow \infty~\forall i \in \mathcal{N}$, $\epsilon > 0$, $V_2(\mathbf{x}_1^*,\mathbf{s}_{d,1}^*)\leftarrow 0$, and $\mathbf{\pi}_1 \leftarrow 0$ 
\WHILE{$\sum_{i \in \mathcal{N}}\sum_{\phi \in \Phi_i}\left(h^{+,\phi}_i + h^{-,\phi}_i \right) > \epsilon$} \label{alg:termination}
\STATE Append $V_2(\mathbf{x}_k^*,\mathbf{s}_{d,k}^*) + \mathbf{\pi}_k^{\intercal} A(\mathbf{x}-\mathbf{x}_k^*) \leqslant \theta$ to~\eqref{eq:master} \label{alg:add_cut}
\STATE Solve master problem~\eqref{eq:master} for $\mathbf{x}_k^*$ \label{op1} \label{alg:solve_M}
\STATE  Given $\mathbf{x}_k^*$, find the worst-case uncertainty realization $s_{d,k}^*$ by solving~\eqref{eq:enumerate} for every scenario in $\hat{\mathcal{U}}$ \label{alg:solve_maxmin}
\STATE Given $\mathbf{x}_k^*$ and $s_{d,k}^*$, solve Subproblem~\eqref{eq:subproblem} for $h^{+,\phi},  h^{-,\phi}, o^{+,\phi}, o^{-,\phi}$ \label{alg:solve_S}
\STATE Set $k \leftarrow k + 1$
\ENDWHILE
\STATE Output: $\mathbf{x}_{k-1}^*, \ o^{+,\phi}, \ o^{-,\phi}$
\end{algorithmic}
\end{algorithm}
\vspace{-.1cm}
\noindent
\\
\textbf{Convergence guarantee:} 
The finite number of blocks in a distribution network guarantees that Algorithm~\ref{alg:cutting_plane} will converge within a finite number of iterations. In extreme cases, sub-gradient cuts will be added until every block is de-energized to meet the termination criterion, as indicated in Step~\eqref{alg:termination}.

\section{Model-free Real-Time Optimal Power Flow}\label{sec:rt_opf}
In this section, we consider real-time optimization of the power injections of DGs within the DNMGs. The goal is to drive the slack bus power injections to follow an active power reference signal $P_0^\bullet$ while keeping voltages within safe bounds as loads change. This is achieved by formulating an OPF problem and using a model-free primal-dual method, adapted for compatibility with DNMGs from~\cite{chen2020model}, to solve it. This approach assumes that generation at any slack bus and all bus voltages can be measured.

Consider a distribution network composed of a set of CCs $\mathcal{E}^\mathrm{CC}$, which represent DNMGs. CC $m$ contains a set $\Phi_m$ of phases and a set $\mathcal{G}_m$ of controllable DGs. 
Let $\mathbf{s}_m$ be a vector of all decision variables corresponding to CC $m$, where each controllable DG $g\in \mathcal{G}_m$ has active and reactive power $s_{g} = \{p_{g}, q_{g}\}$. Let $\Bell_m$ be a vector of the uncontrollable active and reactive power loads in CC $m$. Let $v_m(\textbf{s}_m,\Bell_m)$ be a vector of the voltage magnitudes at every PQ node in CC $m$. Let $P^\mathrm{sl}_{m}(\textbf{s}_m,\Bell_m)$ be the power injection at the slack bus in CC $m$. This work assumes that each CC has one slack bus—either the substation or the bus connected to a GFM-DG—while all other buses remain PQ buses. 

The OPF problem we aim to solve is
\begin{align} \label{eq:control_opt}
    \min & \sum_{m\in \mathcal{E}^\mathrm{CC}}\sum_{\phi\in \Phi_m}\left(P^{sl,\phi}_m(\textbf{s}_m, \Bell_m) - P^{sl,\bullet,\phi}_m \right)^2 \\
    \begin{split}
    \mathrm{s.t.}~ & \underline{v_i} \leqslant v_i^\phi(\textbf{s},\Bell) \leqslant \overline{v_i}, \\ & \qquad \qquad \forall \phi \in \Phi_i,\forall i\in \mathcal{N}_m,\forall m \in \mathcal{E}^\mathrm{CC}
    \end{split} \\
    & 0 \leqslant p^\phi_{g} \leqslant \overline{P}^{\text{DG}}_g, \quad \forall \phi \in \Phi_g,\forall g\in \mathcal{G}_m, \forall m \in \mathcal{E}^\mathrm{CC} \\
    & 0 \leqslant q^\phi_{g} \leqslant \overline{Q}^{\text{DG}}_g, \quad \forall \phi \in \Phi_g, \forall g\in \mathcal{G}_m, \forall m \in \mathcal{E}^\mathrm{CC}
\end{align}
where $ P^{sl,\bullet,\phi}_m$ is an active-power reference signal for slack bus injections in CC $m$ on phase $\phi$; $\underline{v_i}$ and $\overline{v_i}$ are upper and lower voltage limits, respectively; $\mathcal{N}_m$ is the set of buses in CC $m$; $\Phi_i$ is the set of phases at bus $i$; $\Phi_g$ is the set of phases that generator $g$ is connected to; $\overline{P}^{\text{DG}}_g$ and $\overline{Q}^{\text{DG}}_g$ are active and reactive power limits for DG $g$.

\begin{figure}
    \centering
\includegraphics[width=\linewidth]{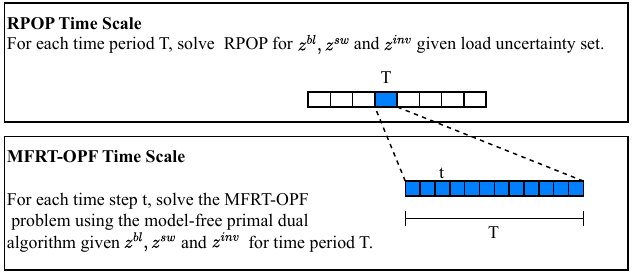}
    \caption{Diagram illustrating the integration of the faster time scale MFRT-OPF algorithm with the RPOP.}
    \label{fig:RPOP-RTOPF}
    \vspace{-0.05in}
\end{figure}

We use the term node to refer to a single phase at a single bus, e.g., a three-phase bus is comprised of three nodes. Define $N^\phi_m$ as the total number of nodes present in CC $m$, i.e., $N^\phi_m := \sum_{i\in\mathcal{N}_m}|\Phi_i|$. Let $\mathbf{\xi^{(k)}} \in \mathbb{R}^{N^\phi_m\times 1}$ be a vector of deterministic exploration signals obtained by sampling 
\begin{equation}
    \xi_g^\phi(t) = \epsilon \cos{(\omega^\phi_g t)}, \quad \forall g \in \mathcal{G}
\end{equation}
where $\epsilon > 0$ is some small value, and $\omega_g^\phi$ is a unique frequency for each generator $g$ and generator phase $\phi$ such that $\omega_i^\phi~\neq~\omega_i^\psi~\neq~\omega_j^\phi~\forall~\phi~\neq~\psi~\in~\Phi_i, \forall i~\neq j \in \mathcal{G}$. Define $f_{m}^{(k)} = \sum_{\phi\in \Phi_m}\left(P^{sl,\phi}_m(\mathbf{x}_m^{(k)}, \Bell_m^{(k)}) - P^{sl,\bullet,\phi}_m  \right)^2$ for each $m\in\mathcal{E}^\mathrm{CC}$ and let
\begin{equation}
    \mathbf{g}^{(k)} = \begin{bmatrix}
        \mathbf{v}(\mathbf{x}^{(k)}, \Bell^{(k)}) - \overline{\mathbf{v}} \\
        \underline{\mathbf{v}} - \mathbf{v}(\mathbf{x}^{(k)}, \Bell^{(k)})
    \end{bmatrix}.
\end{equation}
We introduce a vector $\mathbf{\nu}$ of size $2N^\phi_m$ for each CC. Each element in $\mathbf{\nu}$ is either one or zero. A value of one indicates the node where the associated voltage in \(\mathbf{g}\) is measured is in that CC, while a value of zero indicates the node belongs to a different CC. This vector can be constructed using the switch states and a simple connectivity test.
Let $\rho, \delta > 0$ represent the small regularization parameters of the regularized Lagrangian function, and let \(\alpha > 0\) be a constant step size. Then~\eqref{eq:control_opt} can be solved using Algorithm~\ref{alg:primal-dual}, a modified version of the primal-dual algorithm presented in~\cite{chen2020model}. Figure~\ref{fig:RPOP-RTOPF} illustrates the incorporation of the MFRT-OPF with the RPOP.

\begin{algorithm}[ht]
  \caption{: Model-Free Primal-Dual Algorithm}
  \label{alg:primal-dual}
\begin{algorithmic}[1]
\STATE Initialize $\mathbf{\lambda}_m^{(1)} = \mathbf{0} \in \mathbb{R}^{2N^\phi_m\times 1}$ for each CC $m$ and at each time step k, perform the following steps:
\STATE \textbf{Forward Exploration:} Apply $\mathbf{s}^{(k)}_+ := \mathbf{s}^{(k)} + \mathbf{\xi}^{(k)}$ to the system, and collect measurement $\hat{\mathbf{y}}^{(k)}_+$ of output $\mathbf{y}^{(k)}(\mathbf{s}^{(k)}_+)$.
\STATE \textbf{Backward Exploration:} Apply $\mathbf{s}^{(k)}_- := \mathbf{s}^{(k)} - \mathbf{\xi}^{(k)}$ to the system, and collect measurement $\hat{\mathbf{y}}^{(k)}_-$ of output $\mathbf{y}^{(k)}(\mathbf{s}^{(k)}_-)$.
\STATE \textbf{Approximate Gradient:} Compute the approximate gradient for the primal step 
\begin{equation}
\begin{split}
    \hat{\nabla} \mathbf{\mathcal{L}}_{m}^{(k)}  :=& \nabla_\mathbf{s} \mathbf{f}_m^{(k)}(\mathbf{s}^{(k)}) \\
    &+ \frac{1}{2\epsilon}\mathbf{\xi}^{(k)}\left[ \mathbf{f}_m^{(k)}(\hat{\mathbf{y}}_+^{(k)}) - \mathbf{f}_m^{(k)}(\hat{\mathbf{y}}_-^{(k)}) \right] \\
    &+ \frac{1}{2\epsilon}\mathbf{\xi}^{(k)}(\mathbf{\lambda}_m^{(k)})^\top\biggl[\mathbf{\nu}_m\odot \biggl( \mathbf{g}^{(k)}(\hat{\mathbf{y}}_+^{(k)}) \\ 
    &- \mathbf{g}^{(k)}(\hat{\mathbf{y}}_-^{(k)}) \biggr)\biggr]
\end{split}
\end{equation}
where subscript $m$ denotes values corresponding to CC $m$, $\mathbf{f}_{m}^{(k)}(\hat{\mathbf{y}}^{(k)}) \in \mathbb{R}^{N^{\phi_m} \times 1}$ represents a vector where element $j$ corresponds to $f_{m}^{(k)}$ evaluated at element $j$ of $\hat{\mathbf{y}}^{(k)}$, and $\odot$ represents the Hadamard product. \label{alg:approx_grad}
\STATE \textbf{Approximate Primal Step:} Compute 
\begin{equation}
    \mathbf{s}^{(k+1)} = \text{\it Proj}_{\chi^{(k)}}\{ (1-\alpha \rho)\mathbf{s}^{(k)} - \alpha \sum_{m \in \mathcal{E}^\mathrm{CC}}\hat{\nabla}\mathbf{\mathcal{L}}_m^{(k)} \}  \label{eq:approx_primal}
\end{equation} 
where $\chi$ is a convex set representing operational constraints, and $\text{\it Proj}_{\chi}\{\mathbf{z}\}:=\text{arg} \min_{\mathbf{s}\in\chi} \|\mathbf{z}-\mathbf{s}\|$ denotes projection. \label{alg:approx_primal}
\STATE \textbf{Dual Step:} Compute 
\begin{multline}
    \mathbf{\lambda}_m^{(k+1)} = \text{\it Proj}_{\mathcal{D}^{(k)}}\biggl\{ (1-\alpha \delta)\mathbf{\lambda}_m^{(k)} \\ - 
     \alpha \left(\mathbf{\nu}^{(k)} \odot\mathbf{g}^{(k)}\left(\frac{1}{2}(\hat{\mathbf{y}}_+^{(k)} + \hat{\mathbf{y}}_-^{(k)} )\right)\right)\biggr\} \label{eq:dual_step}
\end{multline} \label{alg:dual_step}
\vspace{-0.15in}
\end{algorithmic}
\end{algorithm}
\vspace{-0.1in}

\section{Case Study: SMART-DS} \label{sec:case_study}
This section presents a case study on a realistic distribution network to demonstrate the effectiveness of both algorithms.

\subsection{Network Details}
\begin{figure}
\centering
\includegraphics[width=0.85\columnwidth]{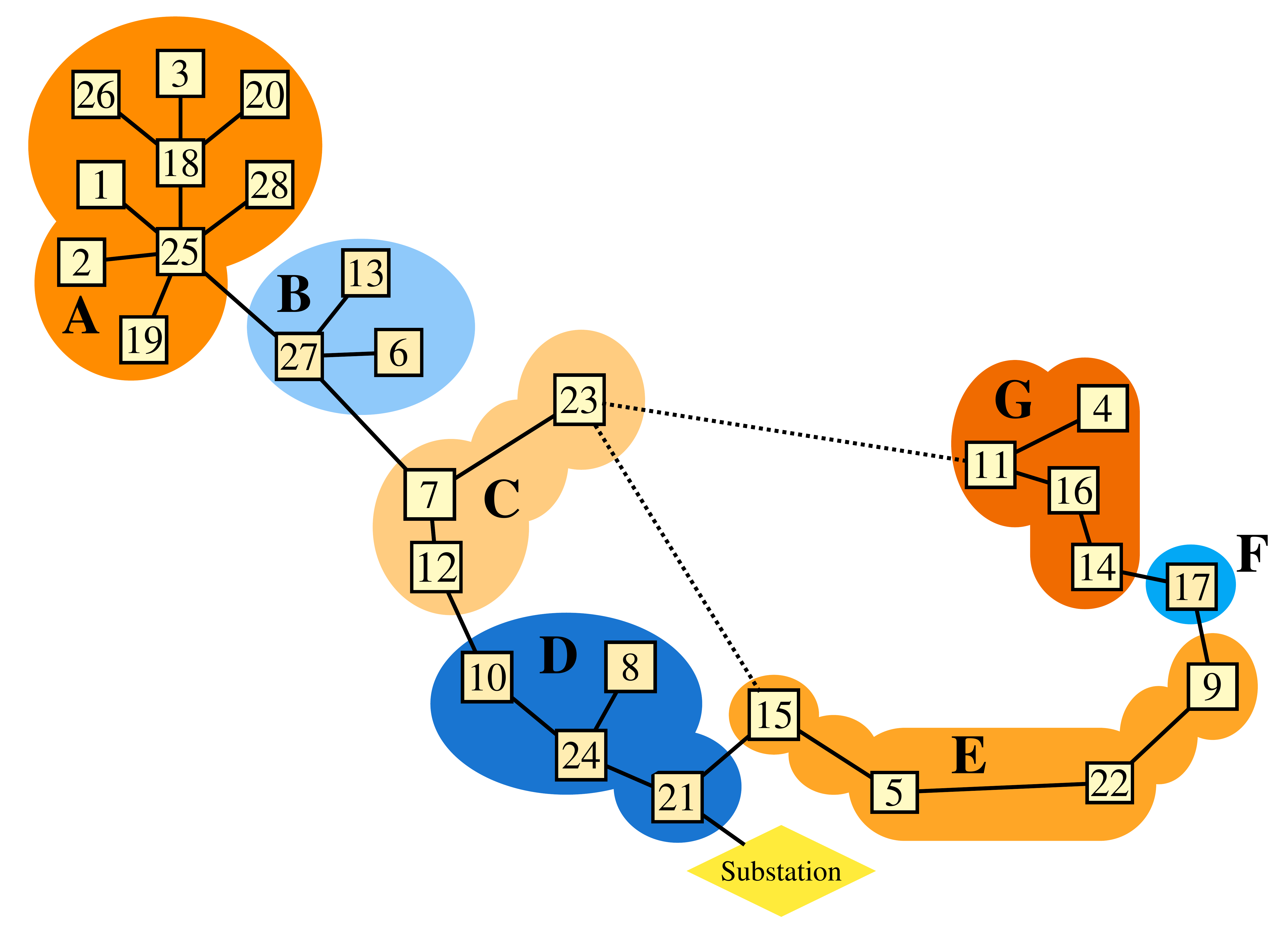}
\caption{Illustration of the load blocks and the switches that connect them. Load clusters A-G are indicated by the blue and orange backgrounds grouping blocks together.}
\label{Block_Map}
\vspace{-.5cm}
\end{figure}

The case study uses a synthetic distribution network based on part of the San Francisco Bay Area from the SMART-DS Synthetic Electrical Network Data OpenDSS models~\cite{SMARTDS}. The network includes $714$ multi-phase buses ($1,565$ total nodes), $652$ lines, and $29$ controllable switches, two of which are redundant lines added to create potential loops. With all controllable switches open, the network forms $28$ blocks, ranging in size from a single single-phase bus to $69$ buses spanning $1$ to $3$ phases. Figure~\ref{Block_Map} shows the block connections. In Fig.~\ref{Block_Map}, blocks are grouped into clusters based on geographic spacing where the total magnitude of active power load is approximately balanced across clusters. When solving the RPOP, it is assumed that loads within each cluster vary together, as described in Section~\ref{subsec:2stage_RPOP}. For example, if a load in Block $6$ is at its maximum, then all loads in Blocks $6$, $13$, and $27$ are also at their maximum.

The network includes $123$ DGs of varying sizes, all assumed to have power injection capabilities like solar photovoltaics. Five DGs are 3-phase, located in Blocks $6$, $17$, $21$, $23$, and $25$, while the rest are single-phase. Blocks $1$, $2$, $3$, $8$, $13$, $19$, $20$, $26$, and $28$ consist only of single-phase buses ($\Phi_l^\text{max} = 1$), while all other blocks contain at least one $3$-phase bus ($\Phi_l^\text{max} = 3$).

\subsection{Implementation Details}
Algorithms~\ref{alg:cutting_plane} and~\ref{alg:primal-dual} were implemented using JuMP v$1.23.2$ in Julia language v$1.11.5$, with PowerModelsONM.jl v3.3.0~\cite{fobes2022optimal} as the foundational framework. Optimization tasks were solved using Gurobi $11.0.2$ for mixed-integer convex problems or IPOPT $3.14.4$ for nonlinear local optimization on a $1.3$ GHz $4$-core Intel Core i$7$ processor with $16$GB memory.

In \eqref{eq:voltLim}-\eqref{eq:PowerBalance}, $\overline{s}_d^{\phi}$, the maximum possible load scenario, is used as the load realization. While any load realization could be applied, using the maximum often led to faster convergence times.

\subsection{RPOP Results} 
\label{subsec:RPOP-results}
\begin{figure*} 
    \centering
  \subfloat[0\% load uncertainty (nominal)\label{fig:nominal_topology}]{%
       \includegraphics[width=0.34\textwidth]{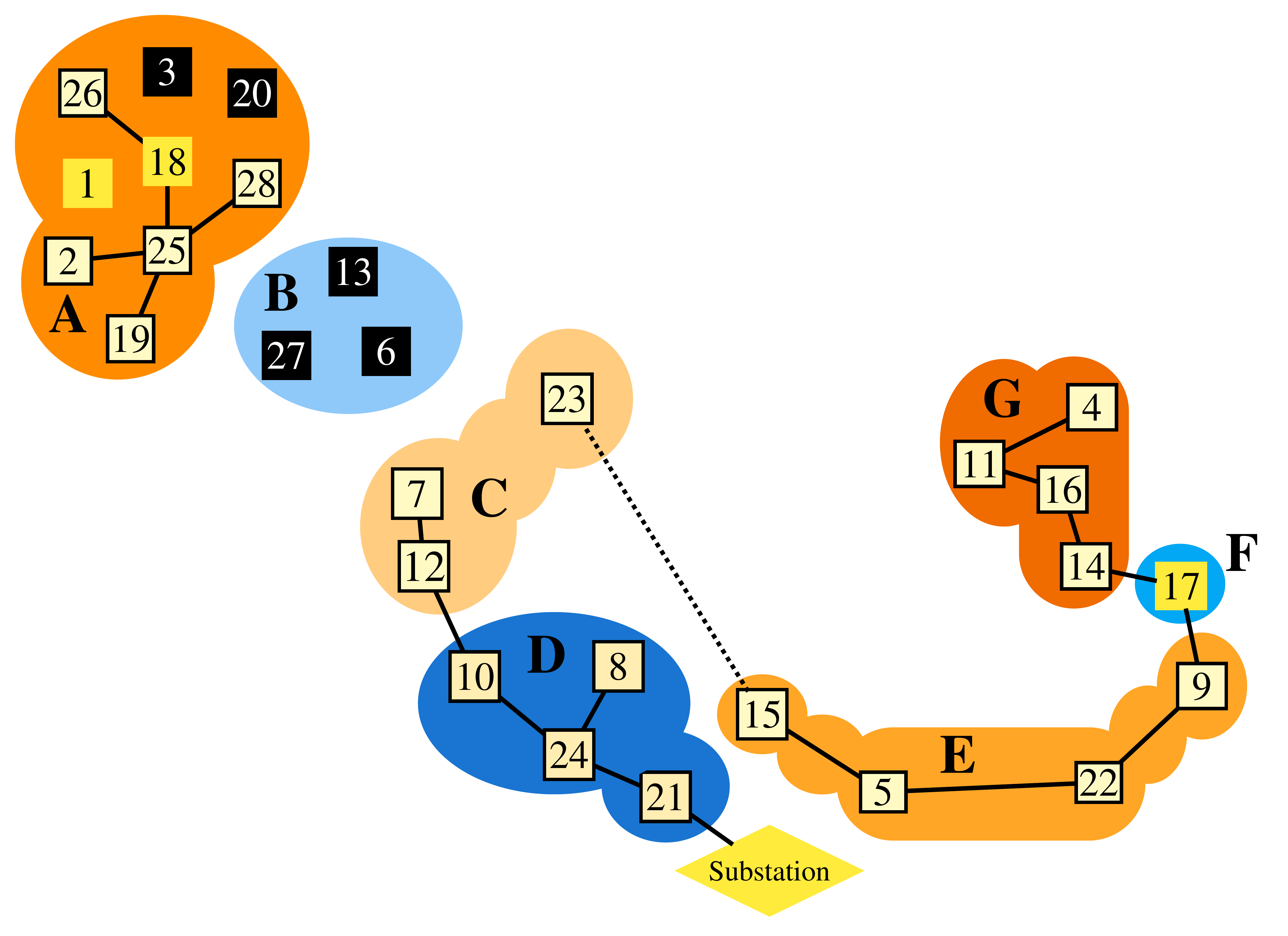}}
    \hfil
  \subfloat[5\% load uncertainty\label{fig:5uc_topology}]{%
        \includegraphics[width=0.34\textwidth]{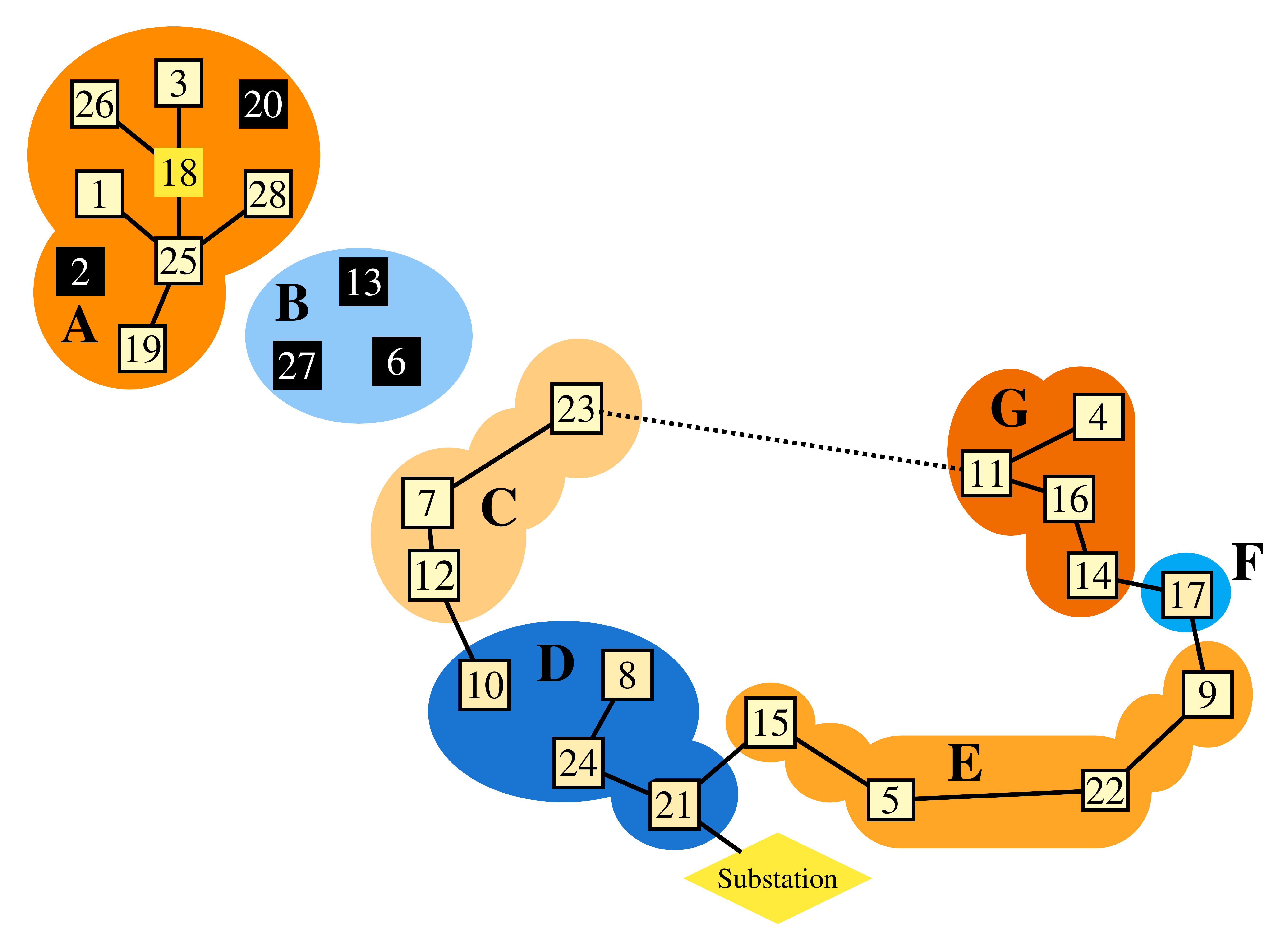}}
    \\
  \subfloat[10\% load uncertainty \label{fig:10uc_topology}]{%
        \includegraphics[width=0.34\textwidth]{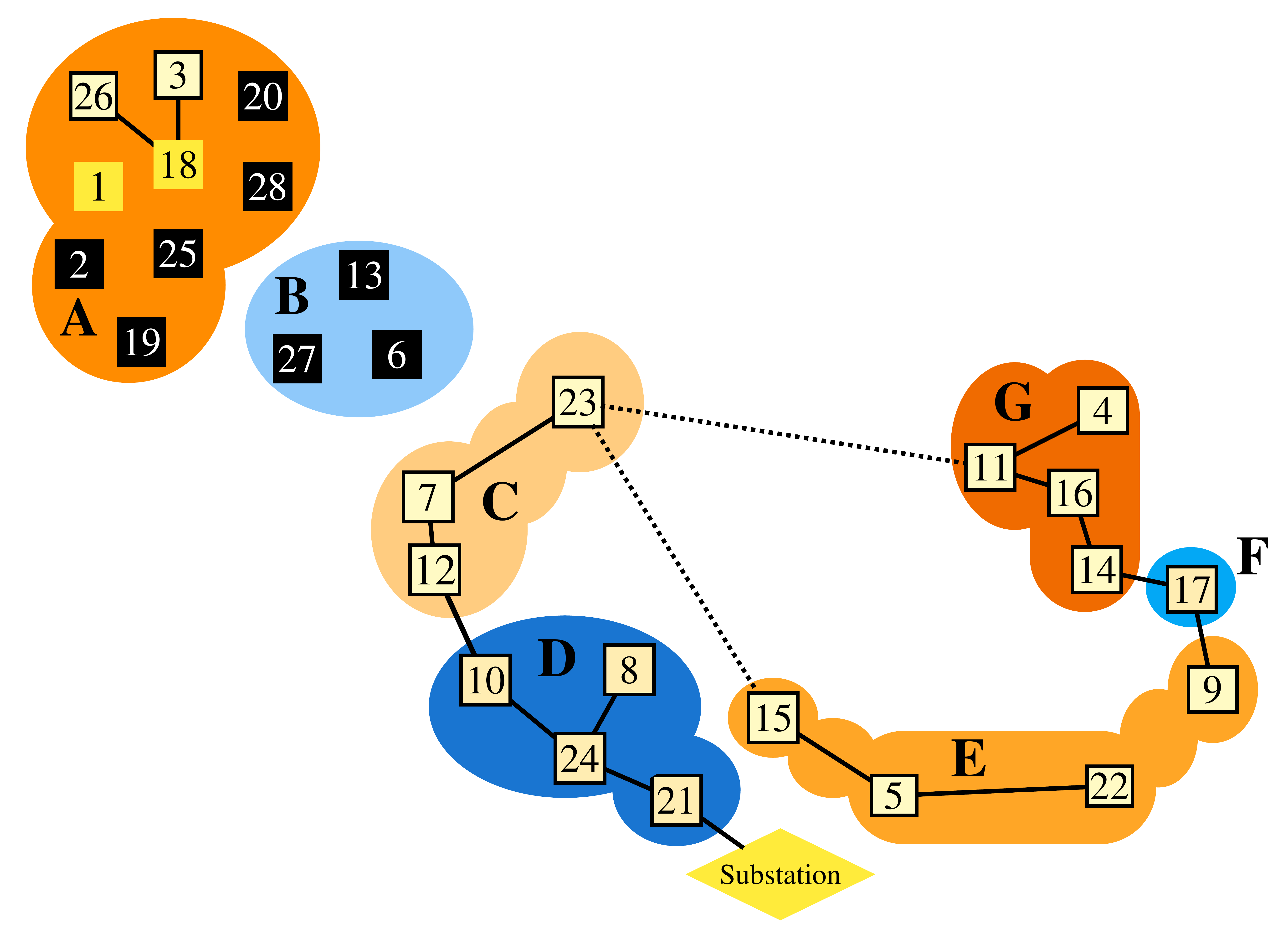}}
    \hfil
  \subfloat[20\% load uncertainty \label{fig:20uc_topology}]{%
        \includegraphics[width=0.34\textwidth]{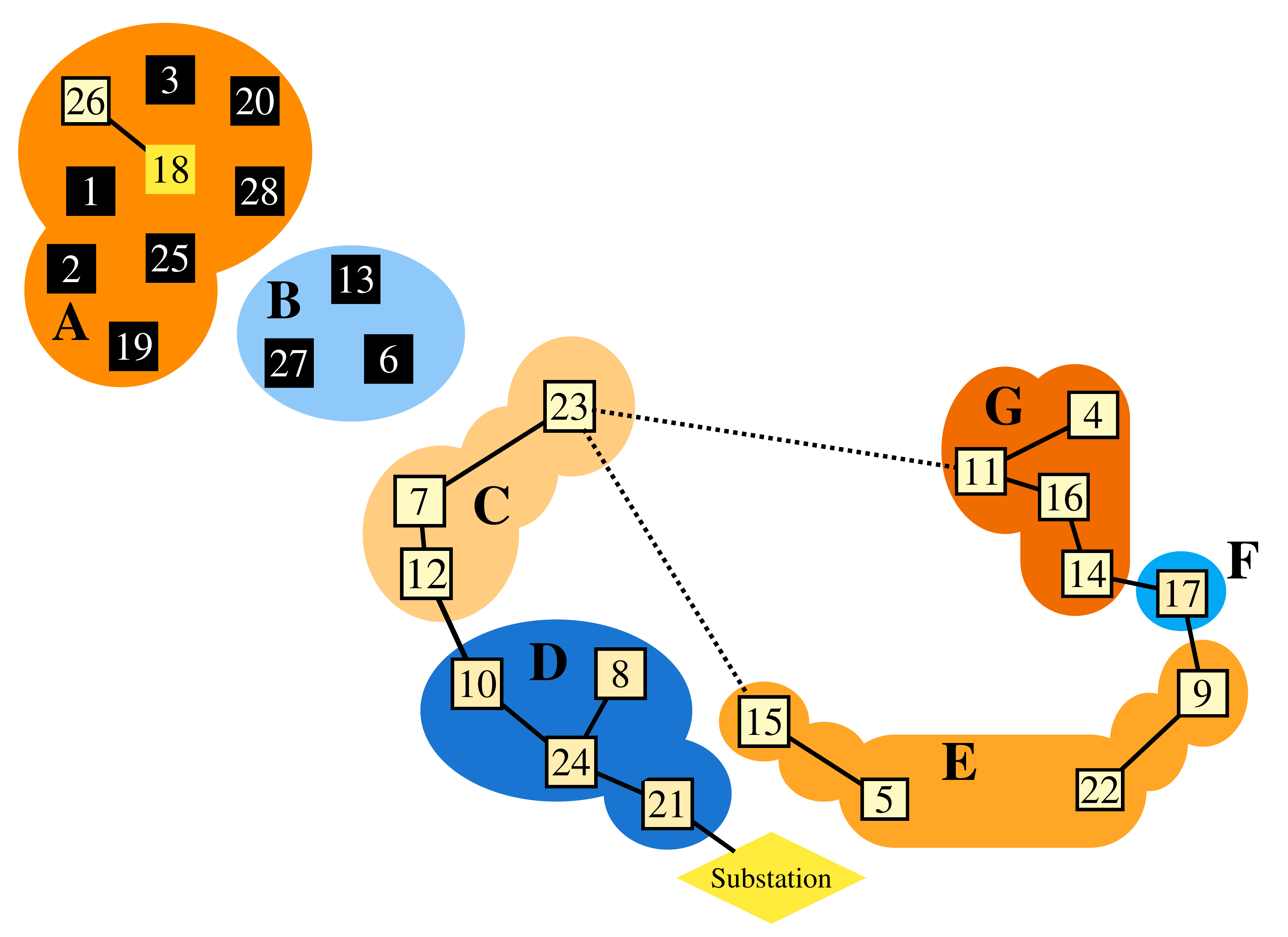}}
  \caption{Optimal partitioning of the SMART-DS network under varying uncertain loads. Black blocks indicate de-energized blocks, yellow blocks without outlines indicate the presence of a GFM-DG or substation voltage source, and black lines represent closed switches. Load clusters A-G are indicated by the blue and orange backgrounds grouping blocks together.}
\label{fig:SmartDS_RPOP_topologies} 
\vspace{-.5cm}
\end{figure*}

The results presented here are for a scenario in which a contingency in Block $27$ requires it to be isolated and disabled. All loads are at their nominal values as published in the SMART-DS dataset. Figure~\ref{fig:nominal_topology} shows the optimal network partitioning without considering load uncertainty. For $5$\%, $10$\%, and $20$\% load uncertainty, the optimal partitioning under this contingency is shown in Figs.~\ref{fig:5uc_topology},~\ref{fig:10uc_topology}, and~\ref{fig:20uc_topology}, respectively. As load uncertainty increases, islanded blocks must be connected to other CCs, and additional blocks may need to be disabled. Run times for the RPOP under each uncertainty level are provided in Table~\ref{tab:run-times}, showing significant reductions compared to the algorithm presented in~\cite{moring2024robust}. This improvement stems from incorporating \textsc{LinDist3Flow} equations with a representative load scenario in the master problem and clustering load uncertainty, reducing the number of cuts (and iterations) needed to approximate the subproblem. The large optimality gap between the new algorithm and the algorithm without linearized power flow in the master problem given in Table~\ref{tab:run-times} illustrates this. For the sake of time, the latter algorithm was forcibly stopped at the start of the first outer iteration after 500 seconds. We omit comparisons with the unclustered case due to its prohibitive computational burden—
solving the unclustered version requires evaluating $2^{|D|}$ scenario subproblems, which is computationally intractable. Furthermore, solutions to the clustered and unclustered problems are not directly comparable: clustering is approximate, and neither version is a guaranteed lower or upper bound of the other.

 \setlength{\tabcolsep}{4pt}
\begin{table}
\renewcommand{\arraystretch}{0.9}
\centering
\caption{RPOP run times}
\label{tab:run-times}
\begin{tabular}{ >{\centering}p{0.20\columnwidth}>
{\centering\arraybackslash}p{0.15\columnwidth}>
{\centering\arraybackslash}p{0.15\columnwidth}>
{\centering\arraybackslash}p{0.15\columnwidth}>
{\centering\arraybackslash}p{0.15\columnwidth}}
\toprule[1.2pt]\midrule[0.3pt]
Load & \multicolumn{2}{c}{Run Time} & \multicolumn{2}{c}{Objective Value} \\
Uncertainty & PF & w/o PF* &  PF & w/o PF \\
(\%) & (s) & (s) & (\$) & (\$)  \\
\midrule
0 & 57.98 &  & &  \\
\midrule
5 & 157.5 & 745.0 & 6.115 & 1.002\\
10 & 215.8 & 505.4 & 8.135 & 1.002\\
20 & 133.9 & 541.8 & 10.17 & 1.002 \\
\midrule[0.3pt]\toprule[1.2pt]
\end{tabular}
\footnotesize{PF = Includes \textsc{LinDist3Flow} in master, w/o PF = No PF in master \\
*these run times reflect the time at which the algorithm was forcibly stopped, i.e., the start of the first outer iteration after 500 seconds.}
\vspace{-0.20in}
\end{table}

 \setlength{\tabcolsep}{2pt}
\begin{table}
\renewcommand{\arraystretch}{0.9}
\centering
\caption{RPOP feasibility}
\label{tab:acpf_robustness}
\begin{tabular}{ >{\centering}p{0.15\columnwidth}>
{\centering\arraybackslash}p{0.20\columnwidth}>
{\centering\arraybackslash}p{0.13\columnwidth}>
{\centering\arraybackslash}p{0.13\columnwidth}>
{\centering\arraybackslash}p{0.13\columnwidth}>
{\centering\arraybackslash}p{0.13\columnwidth}}
\toprule[1.2pt]\midrule[0.3pt]
Load  & Sampled Load & \multicolumn{4}{c}{Feasible Samples} \\
Uncertainty & Uncertainty & \multicolumn{2}{c}{Non-clustered} & \multicolumn{2}{c}{Clustered} \\
& & Lin & Full & Lin & Full \\
(\%) & (\%) & (\%) & (\%) & (\%) & (\%) \\
\midrule
\multirow{3}{*}{0} & 5 & 100 & 100 & 99.8 & 100 \\
& 10 & 100 & 100 & 91.8 & 96.9\\
& 20 & 90.8 & 90.8  & 63.8 & 66.2\\
\midrule
5 & 5 & 100 & 100 & 100 & 100 \\
10 & 10 & 100 & 100 & 100 & 100 \\
20 & 20 & 100 & 100 & 100 & 100 \\
\midrule[0.3pt]\toprule[1.2pt]
\end{tabular}
\footnotesize{Lin = \textsc{LinDist3Flow}, Full = AC-PF}
\vspace{-0.20in}
\end{table}

We empirically assess the robustness of the two-stage RPOP solutions, constrained by linearized \textsc{LinDist3Flow}, against both \textsc{LinDist3Flow} and non-convex AC Power Flow (AC-PF) constraints, while also analyzing the impact of ignoring uncertainty during partitioning. Table~\ref{tab:acpf_robustness} shows the percentage of random load samples from $\mathcal{U}$ that were feasible under \textsc{LinDist3Flow} and AC-PF constraints for different partitioning solutions. These percentages were obtained by solving the RPOP for the assumed load uncertainty level in the first column, fixing the network configuration (switch, inverter, and load block states), and generating $1,000$ random load samples from $\mathcal{U}$ for the uncertainty level in the third column. Each load realization was tested for feasibility with \textsc{LinDist3Flow} and AC-PF equations. Both non-clustered and clustered load uncertainty were sampled, as shown in the $3-4$th and $5-6$th columns, respectively. Samples were generated by randomly selecting values from $[-x, x]$, where $x$ represents the percent uncertainty (e.g., $x=0.2$ for $20$\% uncertainty).
For non-clustered samples, this was applied to each individual load as $s_d = s_d^0(1 + x)$. For clustered samples, the same approach was used, but $x$ was generated per cluster, assigning the same $x$ to all loads in that cluster. Network configuration decisions assume clustered load uncertainty. The AC-PF feasibility values in Table~\ref{tab:acpf_robustness} are approximate lower bounds, since the solver certifies only local infeasibility and each sample was tested from a single warm-start value.

The top three rows in Table~\ref{tab:acpf_robustness}, corresponding to $0$\% load uncertainty during partitioning, show the optimal configuration obtained from the non-robust (deterministic) version of the partitioning problem. These results suggest that accounting for load uncertainty is critical when uncertainty is high but offers limited benefit below $10$\%, especially when loads tend to vary in the same direction, as captured by uncertainty clustering. The bottom three rows, corresponding to partitioning solutions based on the uncertainty levels in the first column, confirm the robustness of the RPOP solution. They also demonstrate that the \textsc{LinDist3Flow} formulation effectively approximates the nonconvex power flow equations.

\subsection{Real-time Control Results} \label{subsec:rt-results}
\begin{figure*} 
    \centering
  \subfloat[Topology for t = $1$ to t = $60$ \label{fig:27out_topology}]{%
       \includegraphics[width=0.30\textwidth]{Figures/20uc_topology_100Load.pdf}}
    \hfil
  \subfloat[Topology for t = $61$ to t = $120$ \label{fig:noContingency_topology}]{%
        \includegraphics[width=0.30\textwidth]{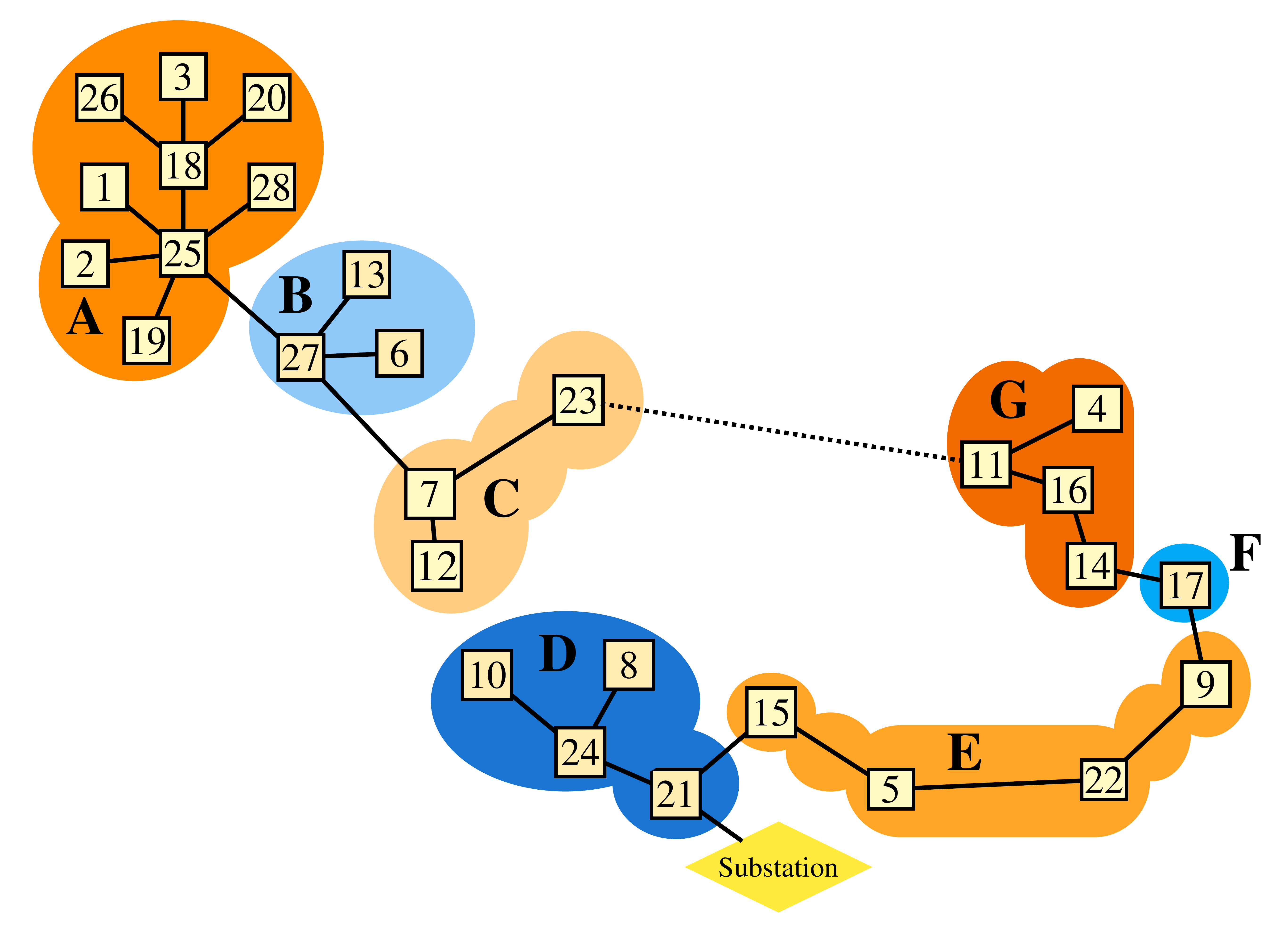}}
    \hfil
  \subfloat[Topology for t = $121$ to t = $180$ \label{fig:SubOut_topology}]{%
        \includegraphics[width=0.30\textwidth]{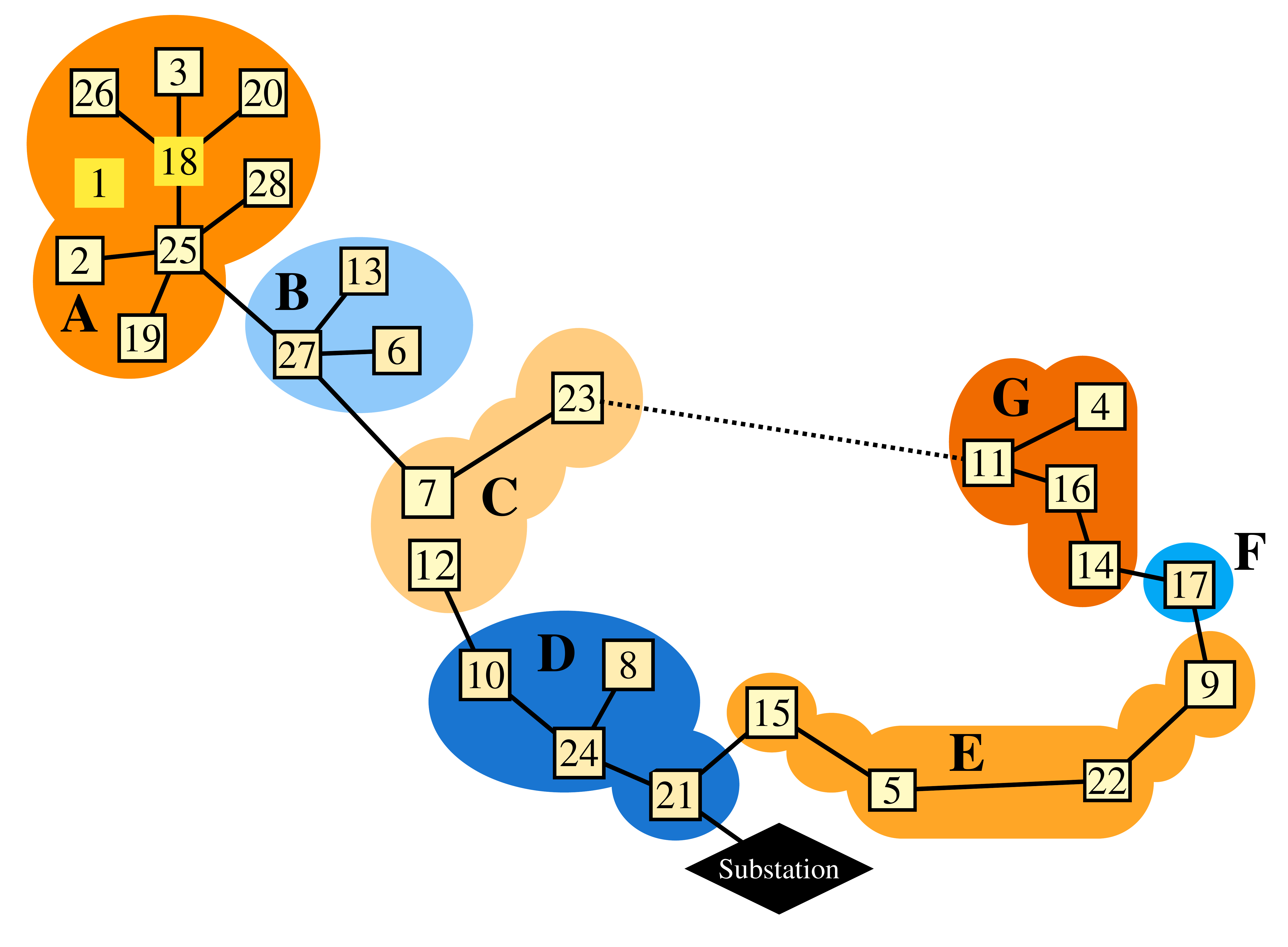}}
  \caption{Optimal partitioning of the SMART-DS network assuming $20$\% load uncertainty. Black blocks indicate de-energized blocks, yellow blocks without outlines indicate the presence of a GFM-DG or substation voltage source, and black lines represent closed switches. Load clusters A-G are indicated by the blue and orange backgrounds grouping blocks together.}
  \label{fig:rt-topologies} 
  \vspace{-.3cm}
\end{figure*}

\begin{figure} 
    \centering
  \subfloat[Total load over time in each CC \label{fig:total_load}]{%
       \includegraphics[width=0.915\columnwidth]{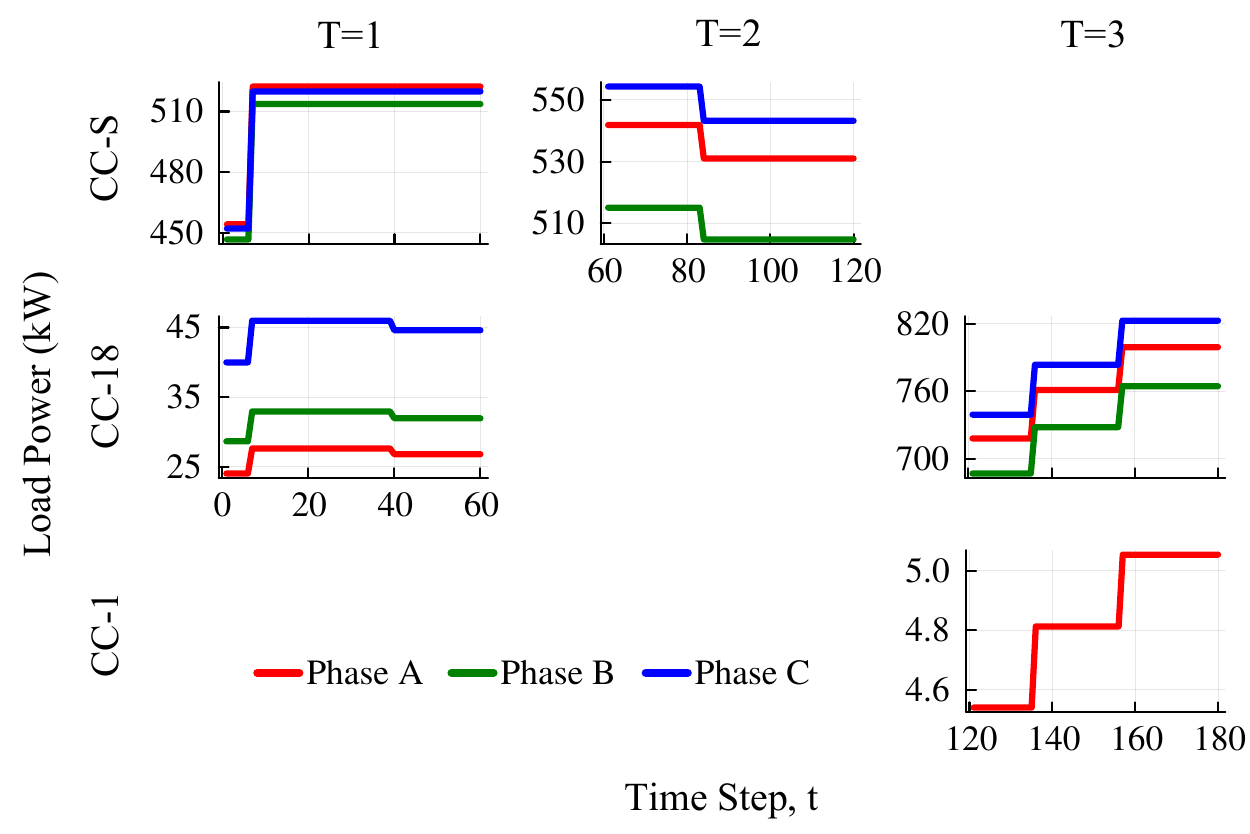}}
\\
  \subfloat[Slack generator power injections over time as network configuration changes \label{fig:slack_injections}]{%
        \includegraphics[width=0.915\columnwidth]{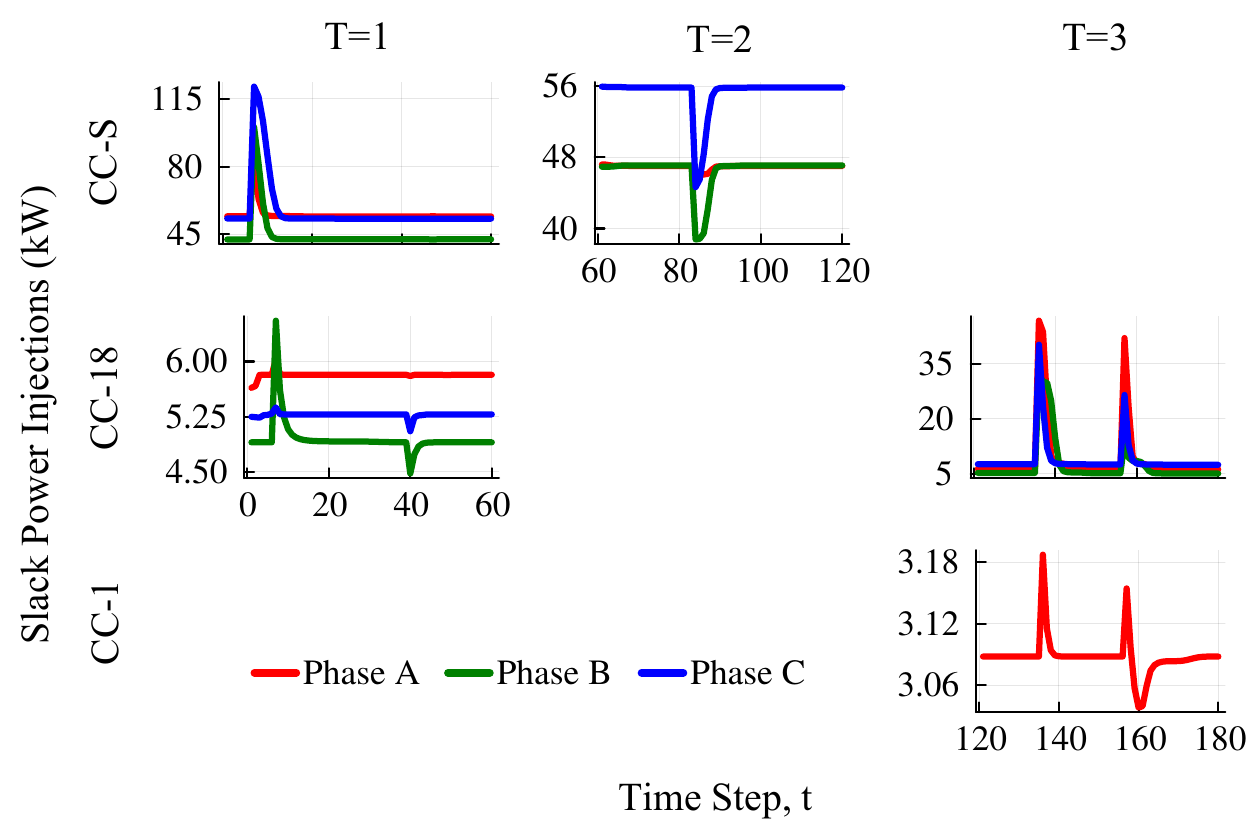}}
\\
  \subfloat[Controlled DG active power injections over time as network configuration changes \label{fig:controlled_power}]{%
        \includegraphics[width=0.915\columnwidth]{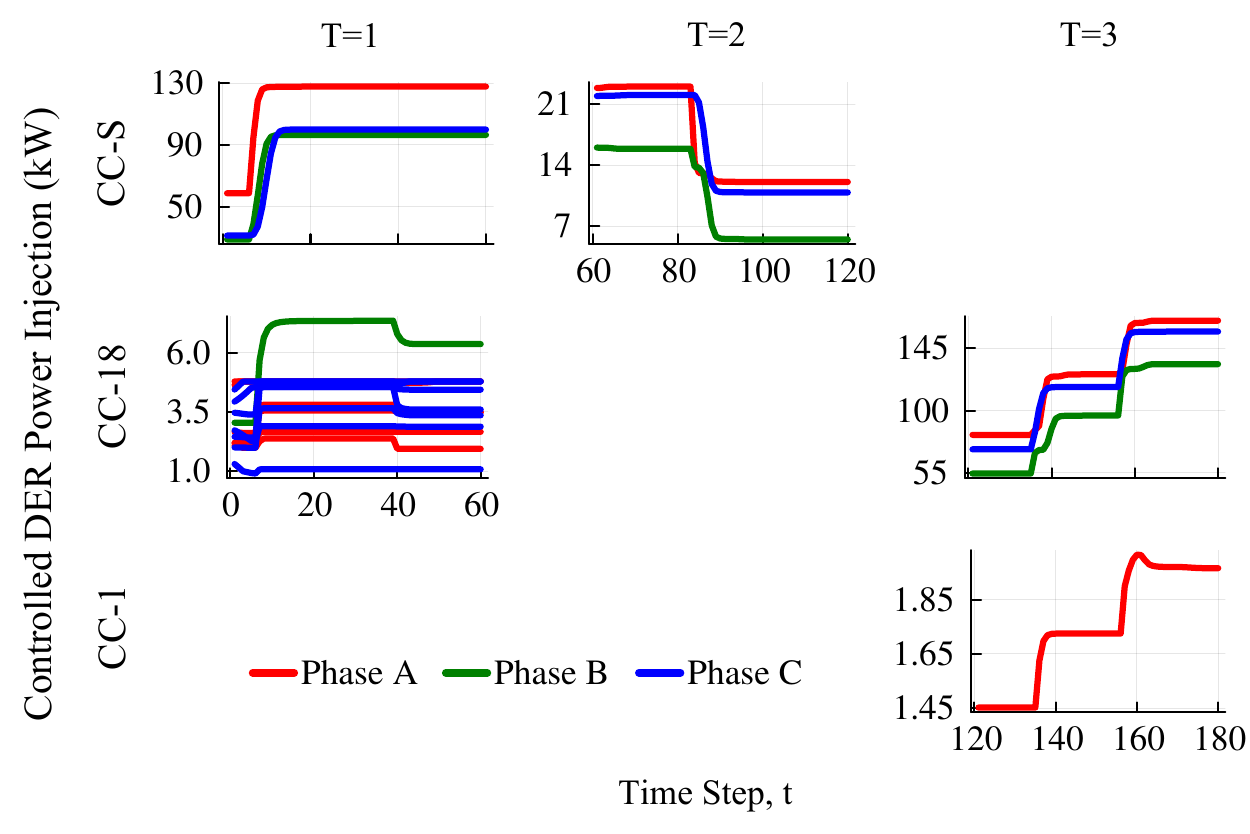}}
  \caption{MFRT-OPF results.}
  \label{fig:rt-results} 
\end{figure}

The results in this section are based on the MFRT-OPF algorithm from Section~\ref{sec:rt_opf}, using the optimal partition from the RPOP under $20$\% load uncertainty across three contingency scenarios. In the first $60$ time steps ($T=1$), Block $27$ is de-energized due to an internal contingency (as in Section~\ref{subsec:RPOP-results}), and all loads start at their nominal values as published in the SMART-DS dataset. From time steps $61$–$120$ ($T=2$), there are no contingencies, and loads start at $90$\% of their nominal values. Finally, from time steps $121$–$180$ ($T=3$), the substation is disconnected, and loads start at $120$\% of their nominal values. The corresponding network partitions for these three scenarios are shown in Figs.~\ref{fig:27out_topology}, \ref{fig:noContingency_topology}, and \ref{fig:SubOut_topology}.

The reference signal for each slack generator, $P_{sl}^\bullet$, is to maintain its initial set-point for each time period. When the network load changes, the slack generator in each CC responds immediately but returns to its original set-point as other controllable generators adjust to address slack injections or voltage violations.

Figure~\ref{fig:total_load} shows the total active power demand in each CC over time. CCs are named after the block hosting the voltage source (e.g., CC-$18$ refers to the CC containing Block $18$ when it includes a GFM-DG, as in $T=1$ and $T=3$). At time step 7, all enabled loads increase by $15$\% from their nominal values. At time step 40, loads in CC-$18$ decrease by $3$\%, while those in CC-S remain unchanged. Additional step changes for enabled loads follow: a $2$\% decrease at time $80$, then $6$\% and $5$\% increases at times $136$ and $157$, respectively. Figure~\ref{fig:slack_injections} shows the slack generators’ responses to load variations. Both generators compensate for the $15$\% increase at time $7$, then return to their nominal set-points. However, when the loads in CC-$18$ are reduced at time step $40$, only the slack generator in Block $18$ responds. Similar response patterns are observed for later load steps under different configurations.

Figure~\ref{fig:controlled_power} illustrates the controlled responses of the DGs to changes in slack generator injections. The set of responsive DGs varies with the network configuration, as the number, size, and load of each CC influence the required capacity across phases. For instance, in $T=1$, CC-$18$ contains only small single-phase DGs (apart from the GFM-DG). To manage potential $20$\% load changes, multiple controllable DGs are needed to restore the GFM-DG to its original set-point.

\section{Conclusions} \label{sec:conclusion}
In this paper, we presented an improved version of the RPOP algorithm presented in~\cite{moring2024robust}. The main improvements were (1) the additional GFM-DG constraints (\eqref{eq:GF_has_max_ph}-\eqref{eq:two-hop}), (2) the inclusion in the master problem of the \textsc{LinDist3Flow} equations for a single load scenario, and (3) the spatial clustering of load uncertainty. Improvement (1) enables the RPOP algorithm to more intelligently select which DG in each CC should be GFM, ensuring that each present phase is connected to a voltage source. Improvements (2) and (3) led to significant improvements in solve time and the ability to use the algorithm on a large, realistic distribution network. 

In addition to the improvements presented for the RPOP, we also presented a MFRT-OPF algorithm. This algorithm enables quick responses to changes in load between reconfigurations. Additionally, the model-free method is well suited for the DNMG setting, where the network structure is frequently changing and may have unknown or uncertain parameters. 

In future work, we aim to incorporate storage units, offering valuable flexibility, into both the RPOP and MFRT-OPF algorithms. For the RPOP, this will require carefully modeling the non-convex complementarity constraints to prevent simultaneous charging and discharging, while ensuring feasible operations for all uncertainty realizations. For both RPOP and MFRT-OPF, the time-coupled nature of state-of-charge will be a key consideration.

\ifCLASSOPTIONcaptionsoff
  \newpage
\fi


\bibliographystyle{IEEEtran}
\bibliography{Bibliography}

\end{document}